\documentclass[preprint,12pt]{elsarticle}

\usepackage{amssymb}
\usepackage{graphicx}
\usepackage{longtable,tabularx}
\usepackage{multirow}
\usepackage{subfigure}
\usepackage{bm}
\usepackage{amsmath,mathrsfs}
\usepackage{cleveref}
\usepackage{xcolor}
\usepackage[title]{appendix}%
\usepackage[absolute,overlay]{textpos}

\crefname{align}{Eq.}{Eqs.}
\crefname{equation}{Eq.}{Eqs.}
\crefname{figure}{Fig.}{Figs.}
\crefname{table}{Table}{Tables}
\crefname{section}{Section}{Sections}
\crefname{appendix}{Appendix}{Appendices}
\journal{Acta Astronautica}

\begin{document}

\begin{frontmatter}



\title{Pontryagin-Bellman Differential Dynamic Programming for Low-Thrust Trajectory Optimization with Path Constraints} %% Article title

\author{Yanis Sidhoum\corref{cor1}}
\author{Kenshiro Oguri}
\affiliation{organization={School of Aeronautics and Astronautics, Purdue University},%Department and Organization
            city={West Lafayette},
            postcode={47907}, 
            state={Indiana},
            country={USA}}

\cortext[cor1]{Corresponding author. Email: ysidhoum@purdue.edu}

%% Abstract
\begin{abstract}
We present a Differential Dynamic Programming framework that parameterizes the control via Pontryagin's Minimum Principle, for constrained low-thrust trajectory optimization. This approach, dubbed Pontryagin-Bellman Differential Dynamic Programming (PDDP), optimizes the costates using a null-space trust-region method, solving a series of quadratic subproblems derived from first- and second-order sensitivities. Terminal equality constraints are handled via a general augmented Lagrangian method, while continuous-time state-path constraints are enforced using a quadratic penalty approach. The resulting solution method represents a significant improvement over classical indirect methods, which are known to face challenges for state-constrained problems. The flexibility of PDDP in handling various state-path constraints is demonstrated through minimum-radius and eclipse-avoidance trajectories in cislunar space. Finally, a trade study against indirect multiple shooting shows that PDDP achieves higher convergence rates overall, indicating improved robustness to poor initial guesses.\end{abstract}

\begin{graphicalabstract}
\end{graphicalabstract}

\begin{highlights}
\item The method combines indirect optimal control with Differential Dynamic Programming 
\item Analytical indirect control law reduces the number of DDP decision variables
\item The hybrid algorithm handles continuous-time state-path constraints
\item Optimizing costates via DDP improves convergence robustness versus indirect shooting
\item Fuel-optimal problem with minimum-radius and eclipse-avoidance constraint is solved
\end{highlights}

%% Keywords
\begin{keyword}


Optimal control \sep Differential dynamic programming \sep Indirect methods\sep Low-thrust trajectory optimization\sep State-path constraints
\end{keyword}

\end{frontmatter}


%% Introduction 
\begin{textblock}{100}(80,10)
\raggedleft
\footnotesize
\textcolor{red}{\textbf{For citation please see published peer-reviewed version: }\\
\emph{Acta Astronautica}, Vol.~247, page 10---24, 2026\\
 \texttt{https://doi.org/10.1016/j.actaastro.2026.04.067}}
\end{textblock}
\section{Introduction}

Low-thrust propulsion systems offer an attractive option for mission designers, as they enable a higher delivered payload mass compared to high-thrust engine alternatives \cite{Benkhoff2017,RACCA20021323}. However, this advantage comes at the cost of long-duration trajectories with throttling, which present significant challenges for mission designers. Over the past few decades, many authors have devoted their efforts on developing efficient and robust solvers, building on the advances in optimal control theory that originated in the early 1960s. Each solver brings its own set of benefits and limitations, making the optimization of low-thrust trajectories an ongoing area of research.

Numerical optimization methods for continuous optimal control problems are typically divided into three categories: 1) indirect methods, 2) direct methods, and 3) differential dynamic programming (DDP) \footnote{We make this distinction based on the underlying principles of optimality: Pontryagin's Minimum Principle for indirect methods, Karush–Kuhn–Tucker conditions for direct methods, and Bellman’s Principle of Optimality for DDP. We acknowledge, however, that this categorization is not a clear-cut—DDP is sometimes viewed as sharing characteristics with both indirect and direct methods.  Our use of these terms is intended to provide context rather than to enforce a strict taxonomy.}. Indirect methods are based on the necessary conditions of optimality derived from the calculus of variations and Pontryagin's Minimum Principle (PMP) \cite{PONTRYAGIN11141962}, leading to an analytical optimal control law known as the \textit{primer vector control law} \cite{LAWDEN1963}, and a two-point boundary value problem (TPBVP) \cite{HULL2003}. They benefit from lower-dimensional search spaces (fewer than 10 variables for 3D control), however, suffer from high sensitivity to initial guesses, leading to a small convergence radius \cite{BERTRAND2002171,RUSSELL2007}. 
Direct methods discretize the problem and solve it via nonlinear programming, offering robustness at the cost of large-scale optimization (typically from hundreds to thousands of variables) and without guarantee of optimality and/or feasibility of the solution \cite{HARGRAVESs1987,Conway2010}. DDP techniques, based on Bellman's Principle of Optimality \cite{Bellman1957}, rely on successive quadratic approximations of the cost function around a nominal trajectory. The resulting quadratic subproblems are solved sequentially to find the feedback control laws that locally
improve the current solution. Initially developed for unconstrained discrete-time problems by Mayne \cite{Mayne1966ASG}, Jacobson and Mayne \cite{Jacobson1970}, and Gershwin and Jacobson \cite{Gershwin1970}, DDP was later enhanced by Whiffen’s Static/Dynamic Control algorithm, which uses Hessian shifting to ensure convexity and penalty functions to handle constraints \cite{Whiffen2002}. This approach is the current state-of-the-art technology in the Mystic software, a high-fidelity low-thrust trajectory design tool used by the Jet Propulsion Laboratory \cite{Whiffen2006}, notably for NASA's Dawn and Psyche missions. The Hybrid DDP (HDDP) algorithm further advances constraint handling with Augmented Lagrangian and active-set methods,  while enabling parallel computing by decoupling optimization from dynamics via first- and second-order State Transition Tensors (STTs) \cite{Lantoine2012Part1,Lantoine2012Part2}. In subsequent works, HDDP was combined with the Sundman transformation \cite{Aziz2018,AzizJGCD}, or with multiple-shooting techniques to form the Multiple-Shooting DDP (MDDP) algorithm \cite{PELLEGRINIMDDPPart1,PELLEGRINIMDDPPart2}.

Building on the long tradition of DDP-based solvers, this work incorporates PMP into the DDP formulation, resulting in a \textit{Pontryagin-Bellman Differential Dynamic Programming} approach, termed PDDP. The connection between DDP and PMP is well established in the literature. Prior studies have shown that the sensitivities of the performance index with respect to the states are generally equivalent to the co-states \cite{Dreyfus1965}. In fact, both quantities satisfy the same differential equations \cite{Jacobson1970}. This relationship has also been explored in the context of collocation (direct methods) and HDDP, where the initial co-states are approximated from the sensitivities of the cost function at the initial point \cite{Lantoine2012Part2}. This procedure is implemented in NASA's Copernicus, a trajectory optimization software \cite{Ocampo2010}. The contribution of this paper is not to revisit this well-known connection. Instead, we parameterize the control using the primer vector control law—explicitly derived from PMP—and solve the constrained optimal control problem with DDP. Previous works have similarly used costate and PMP to parametrize the thrust, however solved the problem with direct methods \cite{Pierson1994}. As a matter of fact, this approach is implemented in Copernicus \cite{Ocampo2010}. The benefits of combining indirect methods with DDP techniques are threefold: first the nonlinear programming techniques embedded in the HDDP algorithm are adapted to the indirect framework, allowing for the handling of state-path constraints, a feature that pure indirect optimization techniques can hardly accommodate; second, using DDP for the optimization of costate variables provides better robustness against poor initial guess compared to classical indirect methods; third, the analytical control law derived from PMP allows for the generation of multi-spiral low-thrust trajectories with a limited number of stages — and, consequently, fewer optimization variables — compared to previous DDP-based solvers. 

In classical indirect methods, the solution to the TPBVP is typically found using the shooting method. This approach can be implemented in a single-shooting form, which tends to exhibit high sensitivity to the initial guess and poor convergence robustness \cite{JIANG2012245}. Alternatively, it can be implemented in a multiple-shooting fashion, which has been shown to improve convergence robustness by distributing the sensitivities across multiple points along the trajectory \cite{SIDHOUMJGCD2024}. Both approaches, however, face challenges when numerically solving state-path constrained problems. In fact, the necessary conditions derived from indirect methods for state-path constrained problems divide the trajectory into a series of constrained and unconstrained sub-arcs, solving each sub-arc as a separate TPBVP while satisfying transversality conditions \cite{Hartl1995,Bryson1975}. This approach is not easily implemented due to two main drawbacks: 1) it requires \textit{a priori} knowledge of the sub-arc structure (e.g., number of sub-arcs and their durations); and 2) optimality conditions for state-path constraint problems lead to discontinuities at the corners, where the costate and/or the Hamiltonian experience jumps, which complicate numerical implementation. Recent works have attempted to incorporate state-path constraints into pure indirect methods, including 1) smoothing approach across constraint corner discontinuities \cite{Oguri2022PrimerVector,OguriCDCIndirect} using activation functions, and 2) interior penalty function methods \cite{Nurre2025}. However, the first approach is limited to path constraints that explicitly depend on the control vector, while the latter leads to multi-point boundary value problems (MPBVPs) with hundreds of shooting segments, due to the sensitivity near the constraint boundaries.

In contrast to traditional indirect shooting methods, in PDDP terminal equality constraints are handled using an
augmented Lagrangian approach \cite{Hestenes1969,Powell1969}. To handle path constraints, on the other hand, we propose integrating the
constraint violations over the trajectory and incorporating the resulting cost as a quadratic penalty \cite{Courant1945} in the augmented
objective. This approach directly addresses continuous-time constraints—a known challenge in many existing direct methods—while remaining general, independent of primer vector, and applicable to DDP with any state representation. To ensure the existence of second-order derivatives of the path-constraint penalty cost, we employ smoothing
techniques inspired by the field of machine learning \cite{Calafiore2020}. Finally, PDDP avoids the need for an active-set strategy for control bounds,
used in previous DDP algorithms \cite{Lantoine2012Part1,Aziz2018,PELLEGRINIMDDPPart1}, as they are inherently satisfied through the primer vector control law.

This paper develops the theory of PDDP as follows. \cref{Sec:PMPReview} provides a brief review of system dynamics and PMP. The main contribution—incorporating PMP into DDP—is detailed in \cref{Sec:PDDP,Sec:PDDPIter}. In \cref{Sec:PDDP}, the PDDP framework is formulated by: (1) applying Bellman's Principle of Optimality to the discrete, PMP-based, optimal control problem (\cref{Sec:DiscreteProbFormulation,Sec:BellmanPrincipleofOptimality}); (2) describing DDP-based costate optimization (\cref{Sec:StructureofPDDP}); (3) comparing PDDP structure with classical DDP \cite{Lantoine2012Part1,Aziz2018,PELLEGRINIMDDPPart1}; (4) handling nonsmooth bang-bang control via state-of-the-art smoothing techniques \cite{TAHERI20182470} (\cref{Sec:BangBangSmoothing}); and (5) enforcing constraints, in particular continuous-time state-path constraints using penalty methods and activation functions (\cref{Sec:AugmentedLagrangian,Sec:StatePathConstraint}). \cref{Sec:PDDPIter} presents the PDDP iteration in detail, including stage quadratic expansions (\cref{Sec:StageQuadraticExpansion}), costate feedback laws (\cref{Sec:TrustRegion}), and sensitivity update equations (\cref{Sec:StageUpdateEquations}). The outer loop, updating Lagrange multipliers and penalty parameters, is described in \cref{sec:EndofIter}. Three minimum-fuel transfers with fixed time of flight in the Earth-Moon system demonstrate the method (\cref{Sec:NumericalExamples}): (1) single- and multi-revolution Distant Retrograde Orbit transfers, for a trade study against indirect multiple shooting; (2) a Halo-to-Halo transfer with a minimum-radius constraint; and (3) a 9:2 NRHO to 2:1 retrograde resonant orbit transfer with eclipse-avoidance constraint. Finally, conclusions are given in \cref{Sec:Conclusion}.

\section{Pontryagin-based Optimal Control Problem Formulation}
\label{Sec:PMPReview}
\subsection{Low-Thrust Engine Model}
\label{Sec:LTModel}
In this paper, we assume a constant specific impulse engine as our low-thrust propulsion system. The control input, $\bm{u} = u \bm{\alpha}$, consists of the unit vector of thrust direction $\bm{\alpha} \in \mathbb{R}^{3}$, i.e. $ \left\|\bm{\alpha} \right\|_{2}=1$, where $\left\|\cdot\right\|_{2}$ denotes the Euclidean norm of a vector, and the engine throttle $u \in[0,1]$. Denote the maximum thrust level of our low-thrust engine by $T_{\mathrm{max}}$, the mass of the spacecraft by $m$, the thruster specific impulse by $I_{\mathrm{sp}}$, and the standard acceleration of gravity at sea level by $g_{0}$ ($\approx 9.80665\times 10^{-3} \mathrm{km/s^{2}}$). The thrust acceleration magnitude $a_{T}$ and the mass-flow rate $\dot{m}$ are computed as: 
\begin{align}
    a_{T} = \frac{T_{\mathrm{max}}}{m}, \qquad \dot{m} = -\dot{m}_{\mathrm{max}}u
\end{align}
where $\dot{m}_{\mathrm{max}} = T_{\mathrm{max}}/(I_{\mathrm{sp}}g_{0})$ is the maximum mass-flow rate.
\subsection{Generic System Dynamics}
\label{sec:GenericSystemDyanmics}
Let $\bm{x}  \overset{\underset{\triangle}{}}{=} [\bm{x}^\top_{\mathcal{O}},m]^\top\in \mathbb{R}^{n_{x}}$ be the full state vector of the spacecraft, where $\bm{x}_{\mathcal{O}}$ is the orbital state of the spacecraft. The system equations of motion are expressed in a generic form as follows: 
\begin{align}
\label{eq:GenericSystemDynamics}
    \dot{\bm{x}} = \bm{f}(\bm{x},\bm{u},t) = \begin{bmatrix}
        \bm{f}_{0}(\bm{x}_{\mathcal{O}},t)+a_{T}\bm{B}(\bm{x}_{\mathcal{O}})u\bm{\alpha}\\
         -\dot{m}_{\mathrm{max}}u
    \end{bmatrix}
        \end{align}
where $\bm{f}_{0}(\cdot):\mathbb{R}^{n_{x_\mathcal{O}}}\times\mathbb{R}\mapsto \mathbb{R}^{n_{x_{\mathcal{O}}}}$ represents the natural dynamics; $\bm{B}(\cdot)\in\mathbb{R}^{n_{x_{\mathcal{O}}}\times3}$ maps the low-thrust acceleration to the rate of orbital state change. For example in orbital mechanics, the natural dynamics $\bm{f}_{0}(\cdot)$ are often described by the unperturbed two-body problem \cite{JIANG2012245}, perturbed two-body problem \cite{Aziz2018,PELLEGRINIMDDPPart2}, or the three-body problem \cite{BOTTA2023}. \cref{eq:GenericSystemDynamics} describes the low-thrust orbital dynamics in various coordinate systems. For instance, in three-dimensional motion with Cartesian coordinates, $\bm{x}_{\mathcal{O}} = [\bm{r}^\top, \bm{v}^\top]^\top \in \mathbb{R}^{6}$, where $\mathbf{r}$ and $\bm{v}$ denotes the position and velocity in the frame of reference, and $\bm{B} = [0_{3\times3}, I_{3}]$. 
\subsection{Pontryagin's Minimum Principle}
\label{sec:PMP}
 Consider a deterministic optimal orbit transfer, where the spacecraft departs the initial orbital state $\bm{x}_{0}$ and arrives at the final state $\bm{x}_{\mathrm{f}}$ while satisfying boundary conditions at the initial and final time. The generic form of this optimal transfer problem is given by: 
\begin{align}
\label{eq:ContinuousOptimalControl}
\begin{aligned}
\min_{\bm{u},t_{0},t_{\mathrm{f}}} \quad &     J = \int^{t_{\mathrm{f}}}_{t_{0}} \mathcal{L}(\bm{x},\bm{u},t) \mathrm{d} t\\
\textrm{s.t.}\quad& \bm{\dot{x}} = \bm{f}(\bm{x},\bm{u},t)  \\
  &\Psi_{0}(t_{0},\bm{x}_{0})=\bm{x}(t_{\mathrm{0}})-\bm{x}_{0}=\bm{0}, \quad \Psi_{\mathrm{f}}(t_{\mathrm{f}},\bm{x}_{\mathrm{f}})=\bm{x}_{\mathcal{O}}(t_{\mathrm{f}})-\bm{x}_{\mathcal{O},\mathrm{f}}=\bm{0}
\end{aligned}
\end{align}
where $\mathcal{L}$ denotes the Lagrangian cost of the problem, and can take the form $\mathcal{L}=\dot{m}_{\mathrm{max}}u$ for minimum-fuel transfers, while $\mathcal{L}=1$ for minimum-time transfers \cite{SIDHOUMJGCD2024}.

A popular approach to solve this type of optimal control problem is the indirect method, based on PMP, which derives the necessary conditions of optimality from the calculus of variation \cite{PONTRYAGIN11141962}. This approach suggests to introduce adjoint variables, also called costate, $\bm{\lambda} \overset{\underset{\triangle}{}}{=} [\bm{\lambda}_{x_{\mathcal{O}}}^\top,\lambda_{m}]^\top \in \mathbb{R}^{n_{\lambda}}$ (where $n_{\lambda}=n_{x}$) to build the Hamiltonian function as $\mathcal{H}\overset{\underset{\triangle}{}}{=} \bm{\lambda}^\top\bm{f}+\mathcal{L}$. The application of Pontryagin's optimality principle yields the following set of optimality necessary conditions, also referred to as Euler-Lagrange equations, which give: 1) the optimal control $(u^{*},\bm{\alpha}^{*}) = \text{arg}\,\min\limits_{(u,\boldsymbol{\alpha})}\mathcal{H}$; 2) state  $\left(\bm{\dot{x}}=\left(\frac{\partial \mathcal{H}}{\partial \bm{\lambda}}\right)^\top=\bm{f}(\bm{x},\bm{\lambda})\right)$; 3) costate $\left( \bm{\dot{\lambda}}=-\left(\frac{\partial \mathcal{H}}{\partial \bm{x}}\right)^\top=\bm{h}(\bm{x},\bm{\lambda})\right)$ dynamics \cite{PONTRYAGIN11141962}.
Lawden's primer vector control law \cite{LAWDEN1963} is summarized in \cref{eq:LawdenPrimerVector}, where $S(\cdot)$ is a function of the costates and mass variables, and termed as the switching function, while $\bm{p}$ is the primer vector:
\begin{align}
\label{eq:LawdenPrimerVector}
     \bm{\alpha}^{*} = -\frac{\bm{p}}{\|\bm{p}\|_{2}}, \quad u^*= \begin{cases}
 1& \text{if }S(t)>0  \\
 0& \text{if } S(t)<0\\
 \in[0,1]& \text{if } S(t)=0 
\end{cases}
\end{align}
The definition of $S$ and $\bm{p}$, as well as the derivation of the optimal control can be found in Ref.~\citenum{BOTTA2023}. The dynamics of the optimal control problem is completely defined by the state and costate dynamics along with the optimal control in \cref{eq:LawdenPrimerVector}. Therefore the canonical variable $\bm{y} \overset{\underset{\triangle}{}}{=}  [\bm{x}^\top,\bm{\lambda}^\top]^{\top}\in \mathbb{R}^{2n_{x}}$ is governed by the following generic equation of motion:
\begin{align}
\label{Eq:GenericCompleteDynamics}
    \dot{\bm{y}} =F(\bm{y})
\end{align}
where $F(\cdot):\mathbb{R}^{2n_{x}}\mapsto \mathbb{R}^{2n_{x}}$ encompasses the state dynamics under low-thrust propulsion in \cref{eq:GenericSystemDynamics}, together with the costate dynamics and optimal control, i.e., $F(\bm{y})\overset{\underset{\triangle}{}}{=} [\bm{f}^\top(\bm{x},\bm{\lambda}),\bm{h}^\top(\bm{x},\bm{\lambda})]^\top$. Therefore, by applying Pontryagin's Minimum Principle the optimal control policy of a low-thrust orbit transfer is parameterized by the costate variables. In classical indirect methods, the solution method consists of finding the initial costate vector $\bm{\lambda}_{0}$ that together with the initial state, satisfies the terminal constraints using the shooting method, a first-order technique \cite{BERTRAND2002171}. Unlike traditional shooting methods, the PDDP algorithm introduced in this work is aimed at solving the costates history, which leads to the optimal control solution, using DDP techniques. The terminal constraints in \cref{eq:ContinuousOptimalControl} are handled using augmented Lagrangian techniques \cite{Hestenes1969,Powell1969}, while state-path constraints are added into the problem formulation in  \cref{eq:ContinuousOptimalControl} and enforced using penalty methods \cite{Courant1945}. The theoretical framework of PDDP is developed in \cref{Sec:PDDP,Sec:PDDPIter}.

\section{Pontryagin-Bellman Differential Dynamic Programming}
\label{Sec:PDDP}
 \subsection{Discrete Problem Formulation}
 \label{Sec:DiscreteProbFormulation}
The first step towards the development of the PDDP algorithm is to discretize the optimal control problem formulated in \cref{eq:ContinuousOptimalControl} into $N$ \textit{stages} (or nodes), focusing on a single-phase formulation. The discrete trajectory under primer vector control law is therefore described by the sequence $[\bm{y}_{1},\cdots,\bm{y}_{N}]$, where $\bm{y}_{k} = \left[\bm{x}^\top_{k},\bm{\lambda}^\top_{k}\right]^\top \in \mathbb{R}^{2n_{x}}$ are the states and costates at stage $k$. The discrete constrained optimal control is formulated as follows: 
\begin{subequations}
\label{Eq:DiscreteOptimalControl}
\begin{equation}
\label{Eq:DiscreteOptimalControlCostFunction}
\min_{\bm{\lambda}_{1},\dots,\bm{\lambda}_{N}} \quad     J = \sum_{k=1}^{N}\mathcal{L}_{k}(\bm{x}_{k},\bm{\lambda}_{k})+\phi(\bm{x}_{N+1})
\end{equation}
\begin{equation}
\label{Eq:DiscreteOptimalControlDynamic}
    \textrm{Dynamics: }\quad \bm{y}_{k+1} = \varphi_{k}(\bm{y}_{k}) \end{equation}
    \begin{equation}
    \label{Eq:DiscreteOptimalControlPathConstraint}
    \textrm{Path constraints: }\quad g_{k}(\bm{x}_{k},\bm{\lambda}_{k}) \leq 0 \end{equation}
  \begin{equation}
  \label{Eq:DiscreteOptimalControlTerminalConstraint}
     \textrm{Terminal constraints: } \Psi(\bm{x}_{\mathrm{N+1}})=\bm{0}
  \end{equation}
\end{subequations}
where $\varphi_{k}:\mathbb{R}^{2n_{x}} \mapsto \mathbb{R}^{2n_{x}}$ are the transition functions that propagates the canonical variable across each stage, $\mathcal{L}_{k}:\mathbb{R}^{n_{x}}\times\mathbb{R}^{n_{x}}\mapsto \mathbb{R}$ are the stage cost functions, $\phi:\mathbb{R}^{n_{x}}\mapsto\mathbb{R}$ is the terminal cost, $g_{k}:\mathbb{R}^{n_{x}}\times\mathbb{R}^{n_{x}}\mapsto \mathbb{R}$ are the stage path constraints, and $\Psi:\mathbb{R}^{n_{x}}\mapsto \mathbb{R}^{n_{\Psi}}$ are the terminal constraints of dimension $n_{\Psi}$. It is important to highlight that in our problem formulation, we do not need to explicitly include thrust magnitude constraints (unlike in Refs. \citenum{Lantoine2012Part1,Aziz2018,PELLEGRINIMDDPPart1}), as the primer vector control law in \cref{eq:LawdenPrimerVector} inherently ensures its satisfaction. We assume that all the functions are at least twice continuously differentiable, and that their first- and second-order derivatives are available. The transition function $\varphi_{k}:\mathbb{R}^{2n_{x}} \mapsto \mathbb{R}^{2n_{x}}$ can be obtained by integrating the spacecraft equation of motion, and is not restricted to deterministic systems. Following our formulation in \cref{Eq:GenericCompleteDynamics,Eq:DiscreteOptimalControlDynamic} the transition function can be expressed as: 
\begin{align}
    \label{Eq:TransitionFunction}
    \varphi_{k}(\bm{y}_{k}) =\bm{y}_{k}+\int_{t_{k}}^{t_{k+1}}F(\bm{y})\mathrm{d}t 
\end{align}
 \subsection{Bellman's Principle of Optimality}
 \label{Sec:BellmanPrincipleofOptimality}
DDP methods are mainly based on Bellman's Principle of Optimality \cite{Bellman1957}. Instead of considering the total cost over the entire trajectory as provided in \cref{Eq:DiscreteOptimalControlCostFunction}, dynamic programming techniques consider the remaining cost, i.e., from the current point to the final destination, termed as \textit{cost-to-go}. At stage $k$, and noting that current or future low-thrust control do not influence the past, the cost-to-go is expressed as: 
\begin{align}
\label{eq:DefCost2Go}
    J_{k}(\bm{x}_{k},\bm{\lambda}_{k},\dots,\bm{\lambda}_{N}) \overset{\underset{\triangle}{}}{=} \sum_{j=k}^{N}\mathcal{L}_{j}(\bm{x}_{j},\bm{\lambda}_{j})+\phi(\bm{x}_{N+1}) 
\end{align}
According to Bellman's Principle of optimality, if the sequence of costates $\bm{\Lambda}^{*}_{k} = \left[\bm{\lambda}_{k}^{*},\dots,\bm{\lambda}_{N}^{*}\right]$ minimizes $J_{k}(\cdot)$, then the sub-sequence $\bm{\Lambda}^{*}_{k+1} = \left[\bm{\lambda}_{k+1}^{*},\dots,\bm{\lambda}_{N}^{*}\right]$ minimizes $J_{k+1}(\cdot)$. Bellman's Principle of Optimality yields the following recursive equation, which is the foundation of dynamic programming \cite{Bellman1957}: 
\begin{align}  
\label{Eq:RecursiveDynamicProgramming}
J^{*}_{k}(\bm{x}_{k})=\min\limits_{\bm{\lambda}_{k}}\left[\mathcal{L}_{k}(\bm{x}_{k},\bm{\lambda}_{k}) + J^{*}_{k+1}(\bm{x}_{k+1})\right]
\end{align}
 \subsection{Structure of PDDP}
 \label{Sec:StructureofPDDP}
DDP techniques seek to minimize a local quadratic approximation of the cost-to-go $J_{k}$, around the current solution. Within PDDP framework, an initial guess consists of a sequence of reference costate vectors at each stage, $\Bar{\bm{\Lambda}} = [\Bar{\bm{\lambda}}_{1}, \dots,\Bar{\bm{\lambda}}_{N}]$. Then, the reference trajectory $\Bar{\bm{x}}$ and control policy $\Bar{\bm{\lambda}}$ are constructed by using the transition function in \cref{Eq:TransitionFunction}, i.e., $\left[\Bar{\bm{x}}^\top(t_{k+1}),\Bar{\bm{\lambda}}^\top(t_{k+1})\right]^\top = \varphi_{k}([\Bar{\bm{x}}^\top_{k},\Bar{\bm{\lambda}}^\top_{k}]^\top)$, for $k=1,\dots,N$. The costates at each node are iteratively updated by applying $\delta \bm{\lambda}_{k}$ until optimality conditions are satisfied. These updates are computed by solving a sequence of quadratic subproblems in a \textit{backward sweep}, inspired by Bellman’s Principle of Optimality. Starting from \( k = N \), each \( \delta \bm{\lambda}_k \) is optimized sequentially, taking into account updates from future stages. Once the backward sweep is complete, the costates are updated via $\bar{\bm{\lambda}}_k^+ = \bar{\bm{\lambda}}_k^- + \delta \bm{\lambda}_k$, resulting in intentional discontinuities in the costate trajectory across stages (see \cref{fig:PDDPSchematic}); this may introduce some suboptimality, though it remains consistent with the local quadratic approximations characteristic of DDP methods. Then, the corresponding state and control trajectories are reconstructed in a \textit{forward sweep} by using the transition function in \cref{Eq:TransitionFunction}. The change in state trajectory is recorded as $\delta \bm{x}_k = \bar{\bm{x}}_k^+ - \bar{\bm{x}}_k^-$, and the process iterates until convergence. It is worth noting that while the costates in PDDP are generally discontinuous across stage boundaries due to the discrete update mechanism, they remain continuous within each stage. As mentioned in the next paragraph, our DDP formulation employs a relatively small number of stages, therefore the costate trajectories are continuous over large portions of the transfer. Additionally, the control parameterization via the primer vector control law ensures that the optimality of the resulting control policy is preserved within each subinterval. The primary source of suboptimality arises from the relaxation of transversality conditions, which may end up not satisfying all necessary PMP conditions. However, such a simplification is commonly adopted in direct and hybrid methods to improve convergence \cite{Ocampo2010,Pierson1994}. Moreover, the presence of costate discontinuities in our framework enables the enforcement of state-path constraints in a practical manner, whereas pure indirect approaches would require explicit handling of corner conditions~\cite{Oguri2022PrimerVector}.

\subsection{Comparison of PDDP to classical DDP}
\label{Sec:ComparisonPDDPClassicalDDP}
\begin{figure}[!ht]
\centering
\includegraphics[width=.9\textwidth]{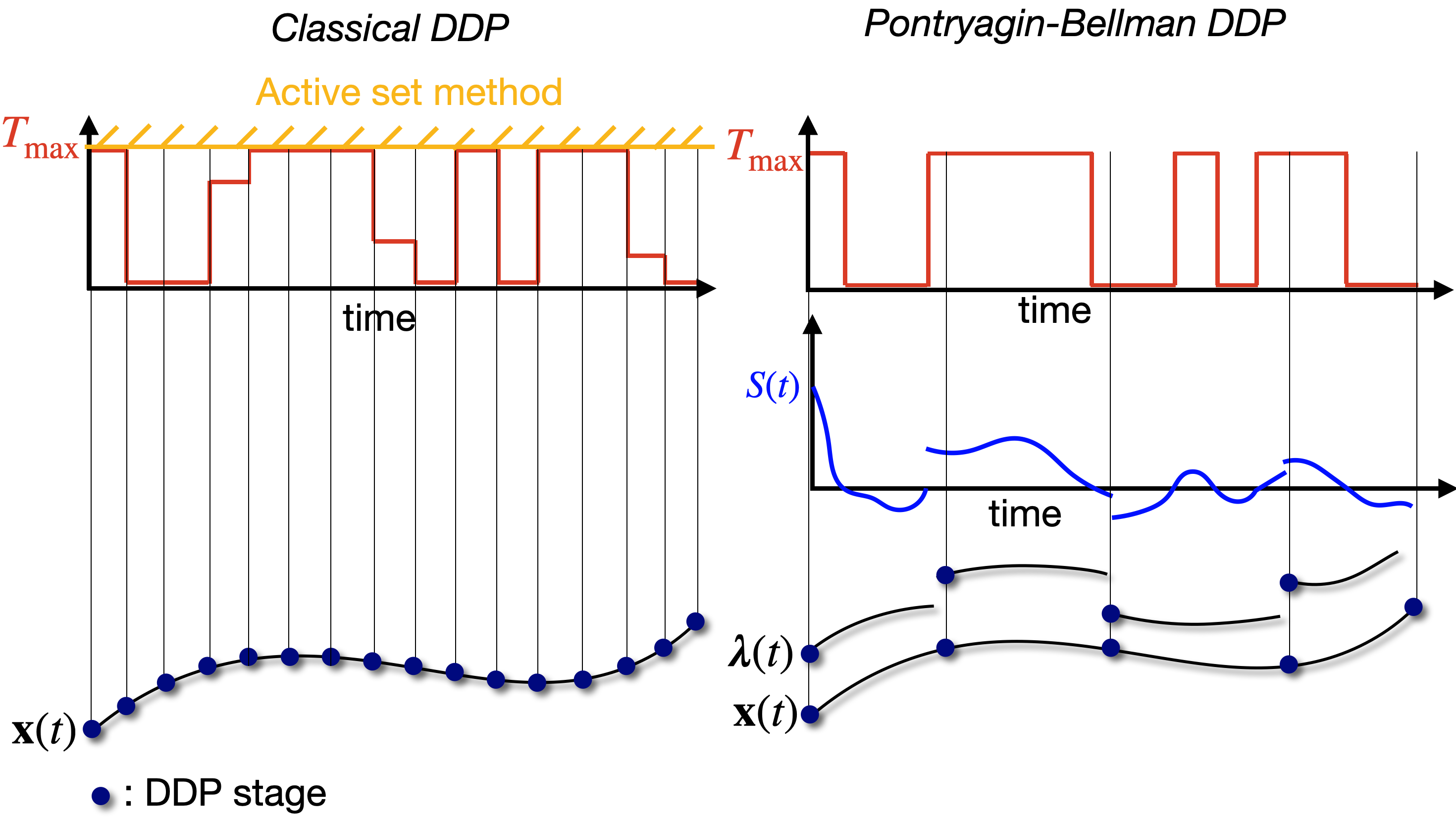}
\caption{Comparison of classical DDP and costate-based control parameterization}
\label{fig:PDDPSchematic}
\end{figure}

Compared to classical DDP, our PDDP framework introduces three enhancements to the general DDP algorithm for long-duration, low-thrust transfers, as illustrated in \cref{fig:PDDPSchematic}. First, in classical DDP, the control $\bm{u}_k$ at each stage is directly optimized under the assumption that it can be approximated by a chosen basis function within the time interval $[t_k, t_{k+1})$, (e.g., piecewise constant as in \cref{fig:PDDPSchematic}, linear, quadratic, sinusoidal, etc). This control approximation necessitates a fine time discretization—often requiring a large number of stages (hundreds to thousands)— to accurately capture changes in control over the trajectory, especially in multi-revolution transfers \cite{Lantoine2012Part2,Aziz2018,PELLEGRINIMDDPPart2}. In contrast, our PDDP approach parameterizes control via costates, which are governed by the Euler-Lagrange equations, allowing bang-bang switches and thrust direction changes to occur naturally between discretization points. As a result, our method can achieve similar or superior accuracy with far fewer discretization stages. Second, classical DDP often yields control profiles that include intermediate thrust levels—values between zero and maximum thrust. While this may lead to near-optimal solutions, it does not enforce the bang-bang structure that characterizes truly optimal low-thrust controls. In contrast, the indirect formulation employed in our PDDP method naturally results in zero or maximum thrust arcs, adhering to the theoretical structure of optimal solutions. The third critical distinction lies in how control constraints are handled. In classical DDP, the control $\bm{u}_k$ is directly optimized at each stage, with thrust magnitude bounds typically enforced using active-set methods \cite{Lantoine2012Part1}. These methods can complicate the algorithm's implementation and occasionally lead to infeasible control updates during the backward sweep. On the other hand, our approach avoids this issue entirely by approximating the control using a smooth hyperbolic tangent function of the costate, as discussed in \cref{Sec:BangBangSmoothing}. This formulation inherently enforces thrust limits by construction, eliminating the need for additional constraint-handling mechanisms and improving numerical robustness.

While our PDDP formulation offers several advantages, it also introduces two main challenges. First, due to the reduced number of discretization points, enforcing state-path constraints only at discrete stages can lead to violations between nodes. To mitigate this, we incorporate a smoothed integral measure of constraint violation directly into the cost function, enabling continuous enforcement as detailed in \cref{Sec:StatePathConstraint}. Second, to retain numerical smoothness in the indirect formulation, we apply a hyperbolic tangent approximation of the throttle, which requires solving a sequence of problems with gradually decreasing smoothing parameters through a homotopy strategy—rather than a single optimization pass.

\subsection{Smooth Approximation of Bang-Bang Control}
\label{Sec:BangBangSmoothing}
The optimal control from PMP in \cref{eq:LawdenPrimerVector} is bang-bang (neglecting singular arcs where $S(t) = 0$ \cite{BERTRAND2002171}), leading to discontinuities in the dynamics \cref{Eq:GenericCompleteDynamics}. Since DDP-based methods require twice continuously differentiable dynamics, we leverage the advances in smoothing techniques for bang-bang control, including homotopy methods \cite{BERTRAND2002171} and hyperbolic tangent smoothing \cite{TAHERI20182470}. As shown in Ref.~\citenum{Sidhoum2023}, hyperbolic tangent smoothing offers infinite differentiability, unlike homotopy, which is not continuously differentiable. For this reason, we use the hyperbolic tangent approximation from Eq. (10) in Ref.~\citenum{TAHERI20182470}, solving the constrained problem over a decreasing sequence of smoothing parameters $\left(\rho_{1,k}\right)_{k}$ until the bang-bang structure emerges. The update strategy for $\rho_1$ within PDDP is detailed in \cref{Sec:Bang-BangSmoothingUpdate}.

 \subsection{Augmented Lagrangian }
 \label{Sec:AugmentedLagrangian}
 The augmented Lagrangian approach, a subset of penalty methods \cite{Hestenes1969,Powell1969}, has been popular to handle terminal equality constraints in DDP frameworks \cite{Lantoine2012Part1,AzizJGCD,PELLEGRINIMDDPPart1}. It transforms a constrained problem into a series of unconstrained ones by adding pure penalty and Lagrange multiplier terms. This allows for much smaller penalty parameters $\sigma$ than pure penalty methods, avoiding the need for $\sigma \rightarrow \infty$ to achieve convergence, thereby reducing the ill-conditioning problems inherent to pure penalty methods \cite{Fletcher2000}. The terminal constraints $\Psi$ in \cref{Eq:DiscreteOptimalControlTerminalConstraint} are relaxed, by adding penalty terms to the terminal cost $\phi$: 
\begin{align}
\label{Eq:TerminalCostAugmentedByTerminalConstraint}
    \hat{\phi}(\bm{x}_{N+1},\bm{\nu}) = \phi(\bm{x}_{N+1})+\bm{\nu}^\top\Psi(\bm{x}_{N+1})+\sigma\|\Psi(\bm{x}_{N+1})\|_{2}^{2}
\end{align}
where $\bm{\nu}\in\mathbb{R}^{n_{\Psi}}$ are the Lagrange multipliers, $\sigma\in\mathbb{R}>0$ is the penalty parameter, and $\hat{\phi}(\cdot):\mathbb{R}^{n_{x}}\times\mathbb{R}^{n_{\Psi}}\mapsto\mathbb{R}$. A recent study surveyed various implementations of augmented Lagrangian methods, with a particular focus on comparing different Lagrange multiplier update strategies—specifically, first- versus second-order updates and inner- versus outer-loop update schemes \cite{Fanger2025}. In \cref{Sec:LagrangeMultipliersUpdate,Sec:PenaltyParameterUpdate}, we briefly summarize the strategy that we use in PDDP. The last step to convert the constrained optimization problem in \cref{Eq:DiscreteOptimalControl} into an unconstrained optimization problem is to relax the state-path constraints in \cref{Eq:DiscreteOptimalControlPathConstraint}. This procedure is detailed in \cref{Sec:StatePathConstraint}. 
 \subsection{Inequality Constraints}
\label{Sec:StatePathConstraint}
As mentioned earlier, we use penalty methods to enforce state-path constraints. We are aware of the numerical limitations of penalty methods in terms of ill-conditioning for large weights \cite{Fletcher2000}. However, as outlined below, in PDDP the stage path constraints $g_k$ are conservatively smoothed for second-order differentiability, and since only soft constraints are enforced, moderate penalty weights are sufficient for convergence. The penalty methods are appealing for low-thrust multi-spiral trajectories, due to the ease of implementation and null effect on the size of the optimization problem. In the framework of low-thrust trajectory optimization, quadratic penalty \cite{Courant1945} have been demonstrated to constrain the thrust magnitude and launch date. In particular, the authors of Refs.~\citenum{Whiffen2002,Aziz2018} incorporate active stage constraints into the stage cost $\mathcal{L}_{k}$, i.e., $\mathcal{L}_{k} = \beta g_{k}(\bm{x}_{k})^{2}$ when $g_{k}(\bm{x}_{k})>0$, while $\mathcal{L}_{k} = 0$ when $g_{k}(\bm{x}_{k})\leq0$ where $\beta$ is the penalty weight. This strategy is deemed possible due to the tight discretization of classical DDP formulation. In contrast, our PDDP formulations involves a much larger discretization, therefore path constraint violations may occur between stages, even if they are satisfied at the stages (see \cref{Sec:ComparisonPDDPClassicalDDP}).

To address this shortcoming, we formulate the penalty term due to path constraints violation differently than the aforementioned studies. Instead of adding the constraint violation at each stage to the stage cost, the accumulated cost $\int_{t_{0}}^{t_{\mathrm{f}}}\ell_{j}(\bm{x}(t),\bm{\lambda}(t))\mathrm{d}t$ due to the violation of the $j$-th path constraint is accounted for over the whole trajectory, where: 
\begin{align}
\label{eq:PathConstViolationAccumulation}
    \ell_{j}(\bm{x}(t),\bm{\lambda}(t)) = \max\left[0,g_{j}(\bm{x}(t),\bm{\lambda}(t))\right], \quad j=1,\dots,n_{c}
\end{align}
where $g_{j}(\cdot):\mathbb{R}^{n_{x}}\times\mathbb{R}^{n_{x}}\mapsto \mathbb{R}$ is the $j$-th state-path constraint, assumed to be twice continuously differentiable, and $n_{c}$ is the number of path constraints. The reader should refer to \cref{eq:Example1PathConstraint,Eq:Eclipseconstraint} in \cref{Sec:Halo2HaloTransfer,Sec:GatewayTransfer}, respectively, for concrete examples of $g(\cdot)$. For conciseness purposes, we include the accumulated path constraint violation in the state vector, i.e., $\bm{x}$ is augmented as: 
\begin{align}
\label{eq:AugmentedVector}
    \dot{\tilde{\bm{x}}}(t) = \begin{bmatrix}\dot{\bm{x}}(t)\\
        \dot{\bm{s}}(t)
    \end{bmatrix} = \tilde{\bm{f}}(\bm{x},\bm{\lambda}) = \begin{bmatrix}
        \bm{f}(\bm{x}(t),\bm{\lambda}(t))\\
        \bm{\ell}(\bm{x}(t),\bm{\lambda}(t))
    \end{bmatrix}
\end{align}
where $\bm{f}(\cdot):\mathbb{R}^{n_{x}}\times\mathbb{R}^{n_{x}}\mapsto\mathbb{R}^{n_{x}}$ is defined in \cref{eq:GenericSystemDynamics} (note that we dropped the explicit time-dependence), and $\bm{\ell}(\cdot):\mathbb{R}^{n_{x}}\times\mathbb{R}^{n_{x}}\mapsto\mathbb{R}^{n_{c}}$ is built by stacking together the constraints in
\cref{eq:PathConstViolationAccumulation}, i.e., $\bm{{\ell}}(\cdot) = [\ell_{1}(\cdot),\dots,\ell_{n_{c}}(\cdot)]^\top$. Therefore, the path constraint violation over the entire trajectory is now accounted by adding $\|\bm{s}(t_{\mathrm{f}})\|_2$ to the terminal  cost $\hat{\phi}$ in \cref{Eq:TerminalCostAugmentedByTerminalConstraint}. In this formulation, $\bm{\ell}(\cdot)$, can be interpreted as the dynamics equation of the state-path constraint violation. It is important to note that $\bm{\ell}(\cdot)$ is continuous, but not differentiable (even to the first order) due to the $\max(\cdot)$ function. However, DDP-based methods require all dynamics equation, path and terminal constraints to be at least twice continuously differentiable. As consequence, a smooth approximation of \cref{eq:PathConstViolationAccumulation} is necessary. 

\begin{figure}[!ht]
\centering
\includegraphics[width=0.5\textwidth]{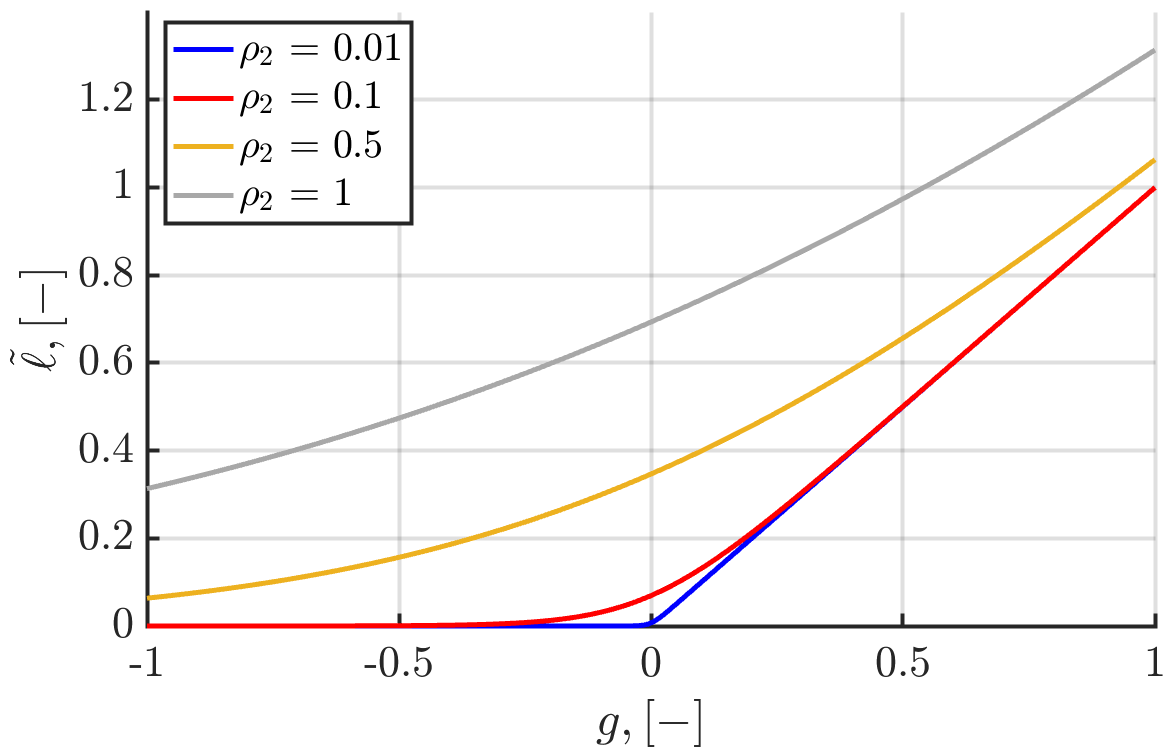}
\caption{Activation function for $\max(\cdot)$ function approximation with various sharpness parameter $\rho_{2}$}
\label{fig:SmoothMaxFunction}
\end{figure}

Activation functions are an intensive topic of research in the field of machine learning, as this class of function allows neural networks to approximate nonlinear, often complex, functions. In this paper, we slightly modify the Log-Sum-Exp (LSE) function \cite{Calafiore2020} to construct a smooth approximation of the $\max(\cdot)$ function in \cref{eq:PathConstViolationAccumulation}. For the sake of conciseness the subscript $j$ is dropped: 
\begin{align}
\label{eq:SmoothApproximationofMax}
        \tilde{\ell}(\bm{x}(t),\bm{\lambda}(t)) = \log{\left(1+\exp{\left(\frac{g(\bm{x}(t),\bm{\lambda}(t))}{\rho_{2}}\right)}\right)}\rho_{2}
\end{align}
where $\rho_{2}\in\mathbb{R}>0$ controls the sharpness of the approximation (smaller values of $\rho_{2}$ lead to sharper approximations). This smooth approximation approach has three important properties that are summarized in Lemma 1. The first property is crucial in the development of the PDDP algorithm. The second property ensures that the approximation preserves satisfaction to path constraints. The third property guarantees the convergence to the local optimum of the original problem. \cref{fig:SmoothMaxFunction} illustrates the behavior of this activation function for different values of the sharpness parameter $\rho_{2}$. Note that, in our implementation $\rho_{2}$ is fixed, although a continuation method over $\rho_{2}$ can be employed.

\textit{Lemma 1:} The smooth approximation in \cref{eq:SmoothApproximationofMax} has the following properties: 
\begin{enumerate}
    \item $\tilde{\ell}(\cdot)$ is twice continuously differentiable across its whole domain 
    \item The approximation in \cref{eq:SmoothApproximationofMax} is conservative, i.e., $\tilde{\ell}(\bm{x}(t),\bm{\lambda}(t))\geq \ell(\bm{x}(t),\bm{\lambda}(t))$
    \item $\tilde{\ell}(\cdot)$ approaches  $\ell(\cdot) $ as $\rho_{2}\to 0^{+}$
\end{enumerate}
\textit{Proof:} See \cref{Sec:ProofLemma1}.

Using \cref{eq:SmoothApproximationofMax} as a smooth approximation of \cref{eq:PathConstViolationAccumulation}, the terminal cost in \cref{Eq:TerminalCostAugmentedByTerminalConstraint} is augmented to form $\tilde{\phi}:\mathbb{R}^{n_{x}+n_{c}}\times\mathbb{R}^{n_{\Psi}}\mapsto \mathbb{R}$ as: 
\begin{align}
\label{Eq:AugmentedTerminalCost}
     \tilde{\phi}(\tilde{\bm{x}}_{N+1},\bm{\nu}) = \phi(\bm{x}_{N+1})+\bm{\nu}^\top\Psi(\bm{x}_{N+1})+\sigma\|\Psi(\bm{x}_{N+1})\|^{2}+\beta \|\bm{s}(t_{N+1})\|_{2}^{2}
\end{align}
where $\beta$ is a penalty weight associated to state-path constraints. The augmented cost-to-go function at stage $k$, $\Tilde{J}_{k}(\cdot)$, of our PDDP method has therefore the following form:
\begin{align}
    \Tilde{J}_{k}(\tilde{\bm{x}}_{k},\bm{\lambda}_{k},\dots,\bm{\lambda}_{N},\bm{\nu})= \sum_{j=k}^{N}\mathcal{L}_{j}(\bm{x}_{j},\bm{\lambda}_{j})+\tilde{\phi}(\tilde{\bm{x}}_{N+1},\bm{\nu})
\end{align}
In the primal-dual framework, the original constrained optimization problem in \cref{Eq:DiscreteOptimalControl} is converted into the following unconstrained minimax problem (see Ref.~\citenum{Bertsekas1982} for more details about primal-dual philosophy): 
\begin{align}
\label{Eq:MinMax}
\begin{split}
&\max_{\bm{\bm{\nu}}} \min_{\bm{\lambda}_{1},\dots,\bm{\lambda}_{N}}\quad     \tilde{J} = \sum_{k=1}^{N}\mathcal{L}_{k}(\bm{x}_{k},\bm{\lambda}_{k})+\tilde{\phi}(\tilde{\bm{x}}_{N+1},\bm{\nu})
    \\
    &\textrm{subject to: }\quad \tilde{\bm{y}}_{k+1} = \tilde{\varphi}_{k}(\tilde{\bm{y}}_{k}) 
    \end{split}
\end{align}
where $\tilde{\bm{y}}\overset{\underset{\triangle}{}}{=}[\tilde{\bm{x}}^\top,\bm{\lambda}^\top] \in \mathbb{R}^{2n_{x}+n_{c}}$, and $\tilde{\varphi}:\mathbb{R}^{2n_{x}+n_{c}} \mapsto \mathbb{R}^{2n_{x}+n_{c}}$ is the augmented counterpart of $\varphi$ defined in \cref{Eq:TransitionFunction}, given by:
\begin{align}
    \label{Eq:AugmentedTransitionFunction}
    \tilde{\varphi}_{k}(\tilde{\bm{y}}_{k}) =\tilde{\bm{y}}_{k}+\int_{t_{k}}^{t_{k+1}}\tilde{F}(\bm{\tilde{y}})\mathrm{d}t 
\end{align}
where $\tilde{F}(\cdot):\mathbb{R}^{2n_{x}+n_{c}}\mapsto \mathbb{R}^{2n_{x}+n_{c}}$ is obtained by stacking together the dynamics of the augmented state $\tilde{\bm{x}}$ and costate $\bm{\lambda}$, i.e., $\tilde{F}(\tilde{\bm{y}})\overset{\underset{\triangle}{}}{=} [\tilde{\bm{f}}^\top(\bm{x},\bm{\lambda}),\bm{h}^\top(\bm{x},\bm{\lambda})]^\top$. In the primal-dual framework, costates $\bm{\lambda}_1, \dots, \bm{\lambda}_N$ are updated in an inner loop to minimize the cost-to-go, while Lagrange multipliers are held fixed and updated in an outer loop to maximize the augmented Lagrangian. Multipliers are initialized to zero, as no intuitive alternative exists, and the state-path constraint variable $\bm{s}$ is also initialized to zero, assuming no initial violation.

\section{PDDP Iterations}
\label{Sec:PDDPIter}
\subsection{Stage-wise Procedure}
\subsubsection{Stage Quadratic Expansion }
\label{Sec:StageQuadraticExpansion}
Denote the augmented state and costate deviations from the nominal solution by $\delta \tilde{\bm{x}}_{k}= \tilde{\bm{x}}_{k}-\Bar{\tilde{\bm{x}}}_{k}\in\mathbb{R}^{n_{x}+n_{c}}$, and $ \delta \bm{\lambda}_{k}= \bm{\lambda}_{k}-\Bar{\bm{\lambda}}_{k}\in\mathbb{R}^{n_{x}}$. Likewise, denote the Lagrange multipliers modification by $\delta \bm{\nu} = \bm{\nu}-\Bar{\bm{\nu}}\in \mathbb{R}^{n_{\Psi}}$. Given that the dynamics equation, terminal, and path constraints are twice continuously differentiable, the augmented cost-to-go function can be expanded to second order in the backward sweep. Assuming that the optimization of the upstream stages has been carried out, the cost-to-go function at stage $k$, $\Tilde{J}_{k}(\Bar{\tilde{\bm{x}}}_{k}+\delta \tilde{\bm{x}}_{k},\Bar{\bm{\lambda}}_{k}+\delta \bm{\lambda}_{k},\Bar{\bm{\nu}}+\delta \bm{\nu})$, is expanded to the second order as: 
\begin{multline}
\label{Eq:QuadraticExpansionCost2Go}
        \delta \tilde{J}_{k} = \tilde{J}^\top_{\tilde{x},k} \delta \tilde{\bm{x}}_{k}+\tilde{J}^\top_{\lambda,k} \delta \bm{\lambda}_{k}+\tilde{J}^\top_{\nu,k} \delta \bm{\nu}+\frac{1}{2}\delta \tilde{\bm{x}}_{k}^\top\tilde{J}_{\tilde{x}\tilde{x},k} \delta \tilde{\bm{x}}_{k}+\frac{1}{2}\delta \bm{\lambda}_{k}^\top\tilde{J}_{\lambda\lambda,k} \delta \bm{\lambda}_{k}\\+\frac{1}{2}\delta \bm{\nu}^\top\tilde{J}_{\nu\nu,k} \delta \bm{\nu}+\delta \tilde{\bm{x}}_{k}^\top\tilde{J}_{\tilde{x}\lambda,k}\delta \bm{\lambda}_{k}+\delta \tilde{\bm{x}}_{k}^\top\tilde{J}_{\tilde{x}\nu,k}\delta \bm{\nu}+
    \delta \bm{\lambda}_{k}^\top\tilde{J}_{\lambda\nu,k}\delta \bm{\nu}+\mathrm{ER}_{k+1}
\end{multline}
where $\delta \tilde{J}_{k} =\Tilde{J}_{k}(\Bar{\tilde{\bm{x}}}_{k}+\delta \tilde{\bm{x}}_{k},\Bar{\bm{\lambda}}_{k}+\delta \bm{\lambda}_{k},\Bar{\bm{\nu}}+\delta \bm{\nu})-\Tilde{J}_{k}(\Bar{\tilde{\bm{x}}}_{k},\Bar{\bm{\lambda}}_{k},\Bar{\bm{\nu}}) $. The subscripts attached to the cost-to-go indicate partial derivatives, for instance $\tilde{J}_{\tilde{x},k}  = \partial\tilde{J}_{k}/\partial \tilde{\bm{x}}_{k} \in \mathbb{R}^{n_{x}+n_{c}} $ and $\tilde{J}_{\tilde{x}\tilde{x},k}  = \partial^{2}\tilde{J}_{k}/\partial \tilde{\bm{x}}_{k}^{2} \in \mathbb{R}^{(n_{x}+n_{c})\times (n_{x}+n_{c})}$. $\mathrm{ER}_{k+1}$ is a constant term which carries the expected quadratic reduction of the objective function from the optimization of upstream stages. In order to find the update in the costate that minimizes the cost-to-go, it is necessary to find the coefficients of this Taylor series expansion. These coefficients are matched to those of the quadratic expansions of $\tilde{J}_{k+1}^{*}(\Bar{\tilde{\bm{x}}}_{k+1}+\delta \tilde{\bm{x}}_{k+1},\Bar{\bm{\nu}}+\delta \bm{\nu})$ and $\mathcal{L}_{k}(\Bar{\tilde{\bm{x}}}_{k}+\delta \tilde{\bm{x}}_{k},\Bar{\bm{\lambda}}_{k}+\delta \bm{\lambda}_{k})$, using the the fundamental recursive equation of dynamic programming $\tilde{J}_{k} = \mathcal{L}_{k}+\tilde{J}_{k+1}^{*}$:  
\begin{align}
\label{eq:QuadraticExpansionStageCost}
&\delta \mathcal{L}_{k}
= \mathcal{L}^\top_{\tilde{x},k} \delta \tilde{\bm{x}}_{k}
 + \mathcal{L}^\top_{\lambda,k} \delta \bm{\lambda}_{k} + \frac{1}{2}\delta \tilde{\bm{x}}_{k}^\top
 \mathcal{L}_{\tilde{x}\tilde{x},k}
 \delta \tilde{\bm{x}}_{k}
 + \frac{1}{2}\delta \bm{\lambda}_{k}^\top
 \mathcal{L}_{\lambda\lambda,k}
 \delta \bm{\lambda}_{k} + \delta \tilde{\bm{x}}_{k}^\top
 \mathcal{L}_{\tilde{x}\lambda,k}
 \delta \bm{\lambda}_{k}\\
\label{eq:QuadraticExpansionOptimCost2Go}
&\delta \tilde{J}^{*}_{k+1}
= \tilde{J}^{*\top}_{\tilde{x},k+1}
 \delta \tilde{\bm{x}}_{k+1}
 + \tilde{J}^{*\top}_{\nu,k+1}
 \delta \bm{\nu}+ \frac{1}{2}\delta \tilde{\bm{x}}_{k+1}^\top
 \tilde{J}^{*}_{\tilde{x}\tilde{x},k+1}
 \delta \tilde{\bm{x}}_{k+1}
 + \frac{1}{2}\delta \bm{\nu}^\top
 \tilde{J}^{*}_{\nu\nu,k+1}
 \delta \bm{\nu}\notag\\& \hspace{8cm}+ \delta \tilde{\bm{x}}_{k+1}^\top
 \tilde{J}^{*}_{\tilde{x}\nu,k+1}
 \delta \bm{\nu}
 + \mathrm{ER}_{k+1}
\end{align}
where $\delta \mathcal{L}_{k} =\mathcal{L}_{k}(\Bar{\tilde{\bm{x}}}_{k}+\delta \tilde{\bm{x}}_{k},\Bar{\bm{\lambda}}_{k}+\delta \bm{\lambda}_{k})-\mathcal{L}_{k}(\Bar{\tilde{\bm{x}}}_{k},\Bar{\bm{\lambda}}_{k}) $, and $\delta \tilde{J}^{*}_{k+1} = \tilde{J}_{k+1}^{*}(\Bar{\tilde{\bm{x}}}_{k+1}+\delta \tilde{\bm{x}}_{k+1},\Bar{\bm{\nu}}+\delta \bm{\nu})-\tilde{J}_{k+1}^{*}(\Bar{\tilde{\bm{x}}}_{k+1},\Bar{\bm{\nu}})$. In order to match the coefficients of \cref{eq:QuadraticExpansionOptimCost2Go} with those of \cref{Eq:QuadraticExpansionCost2Go}, $\delta \tilde{\bm{y}}_{k+1}$ is expressed in terms of $\delta \tilde{\bm{y}}_{k}$ using first- and second-order STTs, $\Phi^{1}_{k} \in \mathbb{R}^{n_{\tilde{y}}\times n_{\tilde{y}}}$ and $\Phi^{2}_{k}\in \mathbb{R}^{n_{\tilde{y}}\times n_{\tilde{y}}\times n_{\tilde{y}}}$, respectively: 
\begin{align}
\label{eq:QuadraticExpansionAugmntedCanonicalVector}
    \delta\tilde{\bm{y}}_{k+1} = \Phi^{1}_{k}\delta\tilde{\bm{y}}_{k}+\frac{1}{2}\delta\tilde{\bm{y}}_{k}^\top\bullet_{2}\Phi^{2}_{k}\delta\tilde{\bm{y}}_{k}
\end{align}
The operator $\bullet_{2}$ denotes a tensor product contracted over the second dimension. For further details on tensor products and the computation of STTs see Ref.~\citenum{Pellegrini2016}. The STTs capture the sensitivities of state and costate variables across stages, enabling the matching of partials in \cref{Eq:QuadraticExpansionCost2Go} with those of $\tilde{J}_{k} = \mathcal{L}_{k} + \tilde{J}_{k+1}^{*}$ via \cref{eq:QuadraticExpansionStageCost,eq:QuadraticExpansionOptimCost2Go,eq:QuadraticExpansionAugmntedCanonicalVector}:
\begin{subequations}
\begin{equation}
     \begin{bmatrix}
            \tilde{J}_{\tilde{x},k}\\
            \tilde{J}_{\lambda,k}
        \end{bmatrix}^\top =  \begin{bmatrix}
            \mathcal{L}_{\tilde{x},k}\\
            \mathcal{L}_{\lambda,k}
        \end{bmatrix}^\top +\begin{bmatrix} \tilde{J}^{*}_{\tilde{x},k+1}\\0_{n_{x}}     
    \end{bmatrix}^\top\Phi^{1}_{k}
    \end{equation}  
\begin{equation}
\begin{aligned}
        \begin{bmatrix}         \tilde{J}_{\tilde{x}\tilde{x},k}&\tilde{J}_{\tilde{x}\lambda,k}\\
           \tilde{J}_{\lambda \tilde{x},k} &\tilde{J}_{\lambda\lambda,k} 
        \end{bmatrix} = \begin{bmatrix}
            \mathcal{L}_{\tilde{x}\tilde{x},k} &\mathcal{L}_{\tilde{x}\lambda,k}\\      \mathcal{L}_{\lambda \tilde{x},k}&\mathcal{L}_{\lambda\lambda,k}
        \end{bmatrix}+\Phi^{1\top}_{k}\begin{bmatrix}
        \tilde{J}^{*}_{\tilde{x}\tilde{x},k+1}&0_{n_{x}+n_{c}\times n_{x}}\\
        0_{n_{x}\times n_{x}+n_{c}}&0_{n_{x}\times n_{x}}     \end{bmatrix}\Phi^{1}_{k}\\+\begin{bmatrix} \tilde{J}^{*}_{\tilde{x},k+1}\\0_{ n_{x}}     
    \end{bmatrix}^\top \bullet_{1} \Phi^{2}_{k}
    \end{aligned}
\end{equation}
\end{subequations}
where the operator $\bullet_{1}$ denotes a tensor product contracted over the first dimension (see Ref.~\citenum{Pellegrini2016} for more details). Finally, the first- and second-order crossed derivatives for the Lagrange multipliers $\nu$ are obtained using the chain rule, noting that $\Phi^{1}_{k}=\frac{\partial \tilde{\bm{y}}_{k+1}}{\partial \tilde{\bm{y}}_{k}}$, and that $\mathcal{L}_{k}$ is independent of $\nu$: 
\begin{align}
\begin{split}
    \tilde{J}_{\nu,k} &= \tilde{J}^{*}_{\nu,k+1}, \quad \tilde{J}_{\nu\nu,k} = \tilde{J}^{*}_{\nu\nu,k+1},\quad [\tilde{J}_{\tilde{x}\nu,k}^\top,\tilde{J}_{\lambda\nu,k}^\top]= [\tilde{J}^{*\top}_{\tilde{x}\nu,k+1},0^\top_{n_{x}\times n_{\nu}}]\Phi^{1}_{k}
\end{split}
\end{align}
In the DDP paradigm, the most computationally intensive task is evaluating first- and second-order STTs \cite{Lantoine2012Part1}. In PDDP, this becomes even more demanding due to the inclusion of costates in the dynamics, effectively doubling the state size. Propagating STTs with the augmented variable $\bm{y}$ requires solving $n_y + n_y^2 + n_y^3$ equations (excluding path constraints). For instance, with $n_y = 14$ (3D case), this results in 2954 equations, compared to 399 in HDDP. To reduce this burden, we exploit tensor symmetry \cite{Pellegrini2016}, lowering the count to $n_y + n_y^2 + n_y^2(n_y + 1)/2$, i.e., 1680 equations for $n_y = 14$.

\subsubsection{Trust Region Method}
\label{Sec:TrustRegion}
The goal of stage-wise quadratic minimization is to compute the costate update $\delta \bm{\lambda}^{*}_{k}$ that minimizes the cost-to-go $\delta \tilde{J}_{k}$ in \cref{Eq:QuadraticExpansionCost2Go}. A naive solution to $\partial \delta \tilde{J}_{k} / \partial \delta \bm{\lambda}_{k} = 0$ may fall outside the valid region of the quadratic expansions or yield non-descent directions if $\tilde{J}_{\lambda\lambda,k}$ is not positive definite. A popular approach to overcome these issues is the trust-region method, which constrains the update within a reliable region by solving the following trust-region quadratic subproblem (TRQP):
    \begin{align}
\label{eq:QuadraticSubproblem}
\text{TRQP}(\tilde{J}_{\lambda,k},\tilde{J}_{\lambda \lambda,k},\Delta):    \min\limits_{\delta\bm{\lambda}_{k}} \tilde{J}_{\lambda,k}\delta\bm{\lambda}_{k}+\frac{1}{2}\delta \bm{\lambda}_{k}^\top \tilde{J}_{\lambda \lambda,k}\delta \bm{\lambda}_{k}, \quad
    \text{s.t.} \quad \left\| D \delta \bm{\lambda}_{k}\right\|_{2} \leq \Delta 
\end{align}
where $\Delta$ is the trust-region radius, and $D$ is a positive-definite scaling matrix that shapes the trust region. Algorithm 7.3.6 in Ref.~\citenum{Conn2000} solves this problem by shifting the Hessian $\tilde{J}_{\lambda \lambda,k}$ by a factor $\gamma^{*}\geq0$ to ensure positive definiteness, i.e.,
$\widetilde{\tilde{J}}_{\lambda \lambda,k}(\gamma^{*}) \triangleq \tilde{J}_{\lambda \lambda,k} + \gamma^{*} D D^\top \succeq 0$. In HDDP convergence was established by replacing $\tilde{J}_{\lambda \lambda,k}$ with its shifted counterpart $\widetilde{\tilde{J}}_{\lambda \lambda,k}(\gamma^{*}) $ in the optimality condition $\partial \delta \tilde{J}_{k} / \partial \delta \bm{\lambda}_{k} = 0$ (see Sec. 5.3.1 in Ref.~\citenum{Lantoine2012Part1} for more details). Applying this wisdom to our PDDP framework, the costate update law is: 
\begin{align}
\label{eq:CostateLaw}
   \delta \bm{\lambda}^{*}_{k}=-\widetilde{\tilde{J}}^{-1}_{\lambda \lambda,k}(\gamma^{*})(\tilde{J}_{\lambda,k} +\tilde{J}_{\lambda \tilde{x},k}\delta \tilde{\bm{x}}_{k}+
   \tilde{J}_{\lambda\nu,k}\delta \bm{\nu})
\end{align}
\subsubsection{Stage Update Equations}
\label{Sec:StageUpdateEquations}
Once the minimization of the quadratic subproblem is completed, as per the costate update law in \cref{eq:CostateLaw}, the sensitivities must be updated accordingly to proceed to the downstream stages. By substituting this costate update control into \cref{Eq:QuadraticExpansionCost2Go}, the expected quadratic reduction and optimized partials at stage $k$ are given by:
   \begin{align}
\begin{split}
    \mathrm{ER}_{k} = \tilde{J}_{\lambda,k}^\top A_{k}+\frac{1}{2}A^{\top}_{k}\tilde{J}_{\lambda\lambda,k}A_{k}+\mathrm{ER}_{k+1}\\
    \tilde{J}^{*}_{\tilde{x},k} = \tilde{J}_{\tilde{x},k} +B^\top_{k}\tilde{J}_{\lambda,k} +A^\top_{k}\tilde{J}_{\lambda\lambda,k}B_{k}+\tilde{J}^\top_{\lambda \tilde{x},k}A_{k}\\
    \tilde{J}^{*}_{\nu,k}=\tilde{J}_{\nu,k} +C^\top_{k}\tilde{J}_{\lambda,k} +A^\top_{k}\tilde{J}_{\lambda\lambda,k} C_{k}+\tilde{J}^\top_{\lambda\nu,k}A_{k}\\
    \tilde{J}^{*}_{\tilde{x}\tilde{x},k} = \tilde{J}_{\tilde{x}\tilde{x},k}+B^\top_{k}\tilde{J}_{\lambda\lambda,k}B_{k}+B^\top_{k}\tilde{J}_{\lambda \tilde{x},k}+\tilde{J}^\top_{\lambda \tilde{x},k}B_{k}\\
    \tilde{J}^*_{\nu \nu,k} = \tilde{J}_{\nu\nu,k} +C^\top_{k}\tilde{J}_{\lambda\lambda,k}C_{k}+C^\top_{k}\tilde{J}_{\lambda\nu,k}+\tilde{J}^\top_{\lambda\nu,k}C_{k}\\
    \tilde{J}^*_{\tilde{x} \nu,k} =\tilde{J}_{\tilde{x}\nu,k}+B^\top_{k}\tilde{J}_{\lambda\lambda,k}C_{k}+B^\top_{k} \tilde{J}_{\lambda\nu,k}+\tilde{J}^\top_{\lambda \tilde{x},k}C_{k}
\end{split}
\end{align}
where $A_{k}\overset{\underset{\triangle}{}}{=}-\widetilde{\tilde{J}}^{-1}_{\lambda \lambda,k}(\gamma^{*})\tilde{J}_{\lambda,k}$ is the feed-forward term, $B_{k}\overset{\underset{\triangle}{}}{=}-\widetilde{\tilde{J}}^{-1}_{\lambda \lambda,k}(\gamma^{*})\tilde{J}_{\lambda \tilde{x},k}$, and $C_{k}\overset{\underset{\triangle}{}}{=}-\widetilde{\tilde{J}}^{-1}_{\lambda \lambda,k}(\gamma^{*})
   \tilde{J}_{\lambda\nu,k}$, are feedback gains in \cref{eq:CostateLaw}. The procedure for solving the quadratic subproblem, as detailed in \cref{Sec:TrustRegion}, along with the updates to the expected quadratic reduction and optimized sensitivities, is recursively repeated in the backward sweep until the first stage. 
   
   At the beginning of the backward sweep at stage $k=N+1$, the optimized sensitivities and the expected quadratic reduction must be initialized. These values are straightforward to compute, as there is no minimization subproblem or stage cost at the final stage $N+1$. Using the general form of terminal cost $\phi$ and constraints $\Psi$ as in \cref{Eq:AugmentedTerminalCost}, we get: 
\begin{align}
\begin{split}
            &\mathrm{ER}_{N+1}= 0\\
         &\tilde{J}^{*}_{\tilde{x},N+1} = \phi_{\tilde{x}}(\bm{x}_{N+1})+\nu^\top \Psi_{\tilde{x}}(x_{N+1})+2\sigma \Psi^\top(\bm{x}_{N+1}) \Psi_{\tilde{x}}(x_{N+1})\\      &\hspace{5.5cm}+2\beta \left[0_{1\times n_{x}},\bm{s}^\top(t_{N+1})\right]\\
        &\tilde{J}^{*}_{\tilde{x}\tilde{x},N+1} = \phi_{\tilde{x}\tilde{x}}(\bm{x}_{N+1})+\nu^\top\bullet_{2}\Psi_{\tilde{x}\tilde{x}}(x_{N+1})+2\sigma \Psi_{\tilde{x}}^\top(x_{N+1})\Psi_{\tilde{x}}(x_{N+1})\\      &\hspace{2.5cm}+2\sigma\Psi^\top(x_{N+1})\bullet_{2}\Psi_{\tilde{x}\tilde{x}}(x_{N+1})+2\beta \begin{bmatrix}
            0_{n_{x}\times n_{x}}&0_{n_{x}\times n_{c}}\\
            0_{n_{c}\times n_{x}} &I_{n_{c}}
         \end{bmatrix}\\
        &\tilde{J}^{*}_{\nu,N+1} = \Psi^\top(x_{N+1}), \quad \tilde{J}^{*}_{\nu \nu,N+1} = 0_{n_{\Psi}\times n_{\Psi}}, \quad\tilde{J}^{*}_{\tilde{x} \nu,N+1}= \Psi_{\tilde{x}}(x_{N+1})
        \end{split}
\end{align}   
\subsection{End of Iteration}
\label{sec:EndofIter}
\subsubsection{Lagrange Multipliers Update}
\label{Sec:LagrangeMultipliersUpdate}
Lagrange multiplier update strategies vary across DDP algorithms: HDDP uses second-order Hessian-based updates \cite{Lantoine2012Part1}, while MDDP adopts first-order steepest ascent \cite{PELLEGRINIMDDPPart1}. A recent study shows that second-order updates in the outer loop yield better robustness to tuning parameters \cite{Fanger2025}. In the outer loop, the multipliers are updated to maximize the augmented Lagrangian by solving TRQP($-\tilde{J}_{\nu,1}, -\tilde{J}_{\nu \nu,1}, \Delta$), yielding:
\begin{align}
\delta \bm{\nu}^{*} = -\widetilde{\Tilde{J}}^{-1}_{\nu \nu,1}(\gamma^{*})J_{\nu,1}=A_{\nu}
\end{align}
where $\widetilde{\Tilde{J}}_{\nu \nu,1}(\gamma^{*})$ is the shifted, negative definite, counterpart of $\Tilde{J}_{\nu \nu,1}$. The expected reduction is updated as:
\begin{align}
\mathrm{ER}{0} = \mathrm{ER}{1} + \tilde{J}^\top_{\nu,1}A_{\nu} + \frac{1}{2}A_{\nu}^{\top}\tilde{J}_{\nu \nu,1}A_{\nu}
\end{align}
\subsubsection{Iteration Acceptance and Trust Region Update}
The quadratic approximations are considered reliable when the ratio  $\chi$ between the actual reduction of the augmented Lagrangian and the predicted quadratic reduction is close to one: 
\begin{align}
    \chi = \frac{\tilde{J}_{\mathrm{new}}-\Bar{\tilde{J}}}{\mathrm{ER}_{0}}
\end{align}
In such cases, the current iterate is accepted, and a larger trust-region size is permitted. However, when the quadratic truncation is deemed unreliable, the iterate is rejected, and the trust-region is reduced to ensure more accurate quadratic approximations. This strategy can be summarized as: 
\begin{align}
\label{eq:TrustRegionRadiusUpdate}
\Delta_{p+1} = \left\{
\begin{aligned}
&\min((1+\kappa)\Delta_{p}, \Delta_{\mathrm{max}}),  \quad &\text{if} \quad \chi \in [1-\epsilon_{1}, 1+\epsilon_{1}], \\
&\Delta_{p},  \quad &\text{if} \quad \chi \in [1-\epsilon_{2}, 1+\epsilon_{2}],\\
&\max((1-\kappa)\Delta_{p}, \Delta_{\mathrm{min}}),  &\text{otherwise}.
\end{aligned}
\right.
\end{align}
where $0<\kappa<1$ is a constant, $p$ is the PDDP iterate, $\Delta_{\mathrm{max}}$ and $\Delta_{\mathrm{min}}$ are the maximum and minimum allowable trust-region size, and $0<\epsilon_{1}<\epsilon_{2}<1$ are user-defined tolerances. 

\subsubsection{Penalty Parameters Update}
 \label{Sec:PenaltyParameterUpdate}
Penalty parameters guide the solution toward feasibility and are updated in the outer loop of the PDDP algorithm, i.e., after accepting an iterate and before the next backward sweep. Their update policy defines both how and when updates occur. Prior works have used different strategies, including min-max formulations \cite{Lantoine2012Part1}, constant penalties \cite{Aziz2018}, and monotonically increasing sequences \cite{PELLEGRINIMDDPPart1}. In PDDP, we adopt a mixed approach: the path constraint penalty $\beta$ remains constant, while the terminal constraint penalty $\sigma$ increases only when constraint violations fail to decrease, continuing until terminal feasibility is achieved:

\begin{align}
\label{eq:PenaltyParameterUpdate}
\sigma_{p+1} = \left\{
\begin{aligned}
&\min(\kappa_{\sigma}\sigma_{p}, \sigma_{\mathrm{max}}),  \quad \text{if} \quad\left(c_{\mathrm{new}}-\Bar{c}\leq \delta_{\sigma}\right) \land \left(c_{\mathrm{new}}\geq\epsilon_{\mathrm{feas}}\right)\\
&\sigma_{p},  \quad \text{otherwise}.
\end{aligned}
\right.
\end{align}
where $c \overset{\underset{\triangle}{}}{=}\|\Psi(\bm{x}_{N+1})\|_{2} $ is the terminal constraint violation, $\kappa_{\sigma}>1$ is a constant, $\sigma_{\mathrm{max}}$ is the maximum allowable penalty parameter, and $\delta_{\sigma}>0$ is a user-defined tolerance that evaluates if the PDDP iterate successfully reduced the terminal constraints. Lastly, $\epsilon_{\mathrm{feas}}$ is a user-defined tolerance for terminal constraints satisfaction. 
\subsubsection{Bang-Bang Smoothing Parameters Update}
\label{Sec:Bang-BangSmoothingUpdate}
As discussed in \cref{Sec:BangBangSmoothing}, the bang-bang optimal control predicted by PMP is smoothed using a hyperbolic tangent approximation \cite{TAHERI20182470}, thereby ensuring that the dynamics are twice continuously differentiable. In PDDP, the optimal control problem defined in \cref{Eq:DiscreteOptimalControl} is solved for a decreasing sequence of smoothing parameters $\rho_{1}$ (see Eq. (10) in Ref.~\citenum{TAHERI20182470} for practical use), i.e., $\rho_{1,0}>\rho_{1,1}>\dots>\rho_{1,\mathrm{min}}$.  The solution from each step serves as the initial guess for the next one, progressively refining the solution until the bang-bang control is captured. The update in the smoothing parameters occurs in the outer loop of PDDP, as follows: 
\begin{align}
\label{eq:SmoothingParameterUpdate}
\rho_{1,p+1} =
\left\{
\begin{aligned}
k_{\rho}\,\rho_{1,p}
&\quad \text{if} \quad
\left(c_{\mathrm{new}} \le \epsilon_{\mathrm{feas}}\right)\land
\left(\mathrm{ER}_{0} \le \epsilon_{\mathrm{opt}}\right)\land \left(
s < \epsilon_{\mathrm{SPC}}\right)\land
\left(\rho_{1,p} > \rho_{1,\mathrm{min}}\right) \\
\rho_{1,p},
&\quad \text{otherwise}.
\end{aligned}
\right.
\end{align}

where $s\overset{\underset{\triangle}{}}{=} \|\bm{s}(t_{N+1})\|$ is the path-constraint violation; $p$ is the PDDP iterate; $0<k_{\rho}<1$ is a multiplying factor; $\epsilon_{\mathrm{opt}}$, $\epsilon_{SPC}$, and $\rho_{1,\mathrm{min}}$ are user-defined tolerances. Specifically, $\epsilon_{\mathrm{opt}}$ is the optimality tolerance, $\epsilon_{SPC}$ the path-constraint violation tolerance, while $\rho_{1,\mathrm{min}}$ represents the threshold at which bang-bang control is achieved. Note that each time the sharpness parameter $\rho_{1}$ is decreased, we reset the penalty parameter $\sigma$ to its initial value $\sigma_0$, since the updated control (due to reduced sharpness) typically increases constraint violations; keeping a large penalty at this point in the iteration would overly bias the solution toward feasibility over optimality.

\subsubsection{Convergence Test}
\label{Sec:ConvergenceTest}
\begin{figure}[!ht]
\centering
\includegraphics[width=1\textwidth]{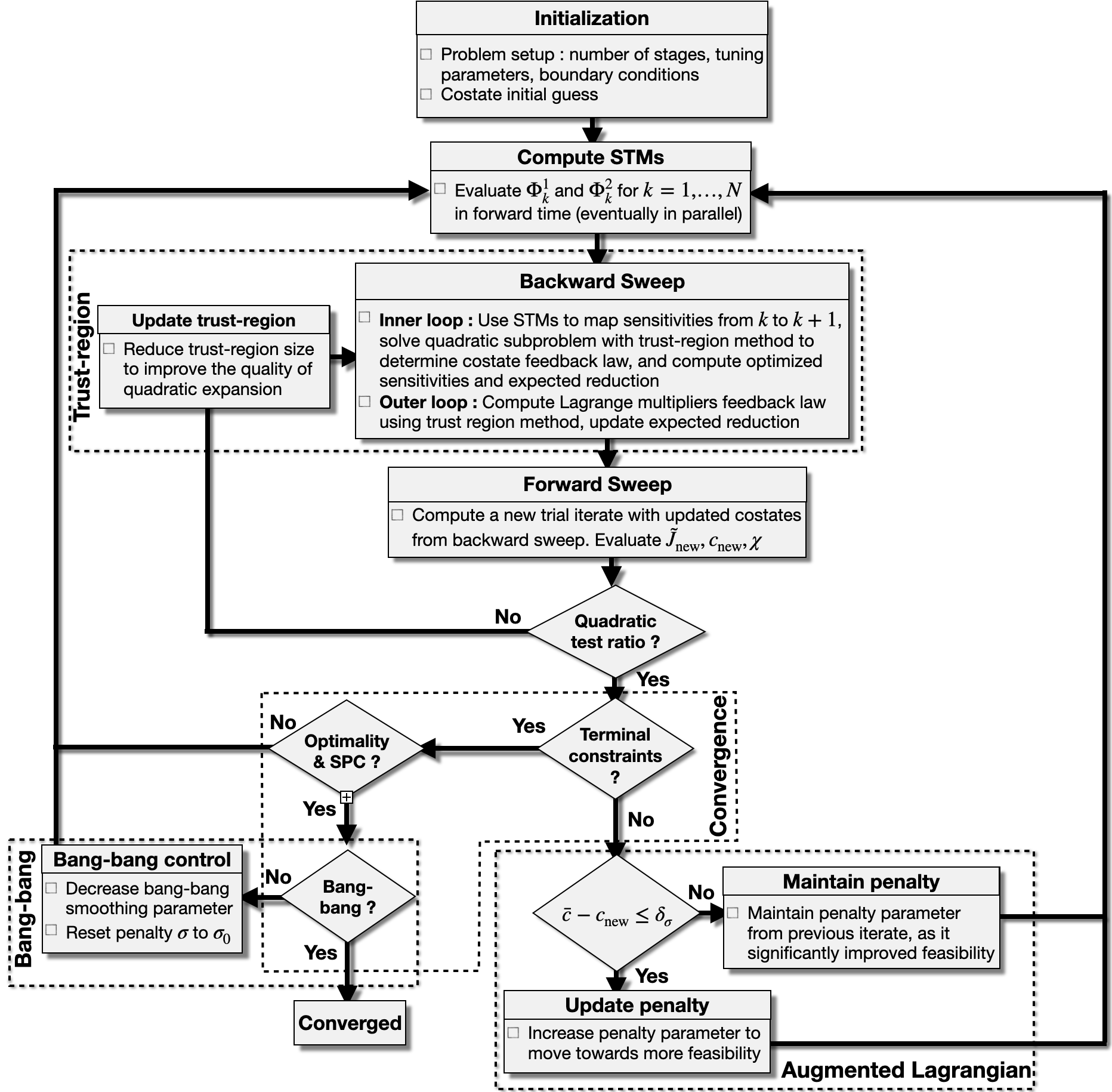}
\caption{PDDP flow chart}
\label{fig:PDDPFlowChart}
\end{figure}

PDDP terminates when the following are satisfied: (i) $\mathrm{ER}0 < \epsilon_{\mathrm{opt}}$ (optimality), (ii) $c < \epsilon_{\mathrm{feas}}$ and $ s < \epsilon_{\mathrm{SPC}}$ (feasibility), (iii) $\rho_1 < \rho_{1,\mathrm{min}}$ (bang-bang control), and (iv) $\Tilde{J}_{\lambda \lambda,k} \succeq 0$, $\Tilde{J}_{\nu \nu,k} \preceq 0$ $\forall k$ (second-order optimality).
In \cref{fig:PDDPFlowChart}, we summarize the flow of the PDDP algorithm, with the different building blocks highlighted. Note that in our implementation, the sharpness parameter $\rho_{2}$ for the approximation in \cref{eq:SmoothApproximationofMax} is fixed over PDDP iterations, however, the user is free to perform a continuation over $\rho_{2}$.

\section{Numerical Examples}
\label{Sec:NumericalExamples}
The theoretical developments on PDDP in \cref{Sec:PDDP,Sec:PDDPIter} are numerically demonstrated by solving three illustrative example problems. The dynamical model used for these examples is the Circular Restricted Three-Body model, with a specific focus on the Earth-Moon system. The characteristics length, time, and mass parameter of the Earth-Moon system are $l^{*}=  385,692.5 $ km, $t^{*}=377084.1527$ s, and $\mu= 0.012150586$, respectively. The dynamics are described using Cartesian coordinates with $\bm{x}_{\mathcal{O}} = [\bm{r}^\top,\bm{v}^\top]^\top \in\mathbb{R}^{6}$ (i.e., $n_{x_{\mathcal{O}}}=6$) as the orbital state of the spacecraft, where $\bm{r} =[
 x, y, z]^\top $ and $\bm{v}=[
 \dot{x}, \dot{y}, \dot{z}]^\top $ denote the position and velocity of the spacecraft in the synodic frame. Therefore the state of the spacecraft is $\bm{x} = [\bm{r}^\top,\bm{v}^\top,m]^\top\in\mathbb{R}^{7}$ (i.e., $n_{x}=7$), while the costate vector is given by $\bm{\lambda} = [\bm{\lambda}_{r}^\top,\bm{\lambda}_{v}^\top,\bm{\lambda}_{m}]^\top\in\mathbb{R}^{7}$ (i.e., $\bm{\lambda}_{x_{\mathcal{O}}}=[\bm{\lambda}_{r}^\top,\bm{\lambda}_{v}^\top]^\top$). The mapping of low-thrust acceleration to the rate of orbital state change is given by $\bm{B}=[0_{3\times3},I_{3}]$. In this framework, the explicit form of the canonical variable dynamics $F$ in \cref{Eq:GenericCompleteDynamics} can be found in Eq. (15) of Ref.~\citenum{BOTTA2023}.
 
In the scenarios studied in the following, we aim at minimizing the fuel consumption, therefore the cost in \cref{Eq:DiscreteOptimalControlCostFunction} does not include local stage costs, i.e., $\forall k$ $ \mathcal{L}_{k}=0$, while the terminal cost $\phi(\bm{x}_{N+1})$, constraints $\Psi(\bm{x}_{N+1})$, and their partials are: 
\begin{align}
\begin{split}
    &\phi(\bm{x}_{N+1})=-m(t_{\mathrm{f}}), \quad \phi_{\tilde{x}}(\bm{x}_{N+1}) = \left[0_{1\times 6},-1,0\right], \quad \phi_{\tilde{x}\tilde{x}}(\bm{x}_{N+1}) = 0_{8\times8},\\
    &\Psi(\bm{x}_{N+1}) = [\bm{r}^\top(t_{\mathrm{f}})-\bm{r}_{\mathrm{f}}^\top,\bm{v}^\top(t_{\mathrm{f}})-\bm{v}^\top_{\mathrm{f}}]^\top, \Psi_{\tilde{x}}(x_{N+1})=\left[
         I_{6},0_{6\times 2}\right], \\&\Psi_{\tilde{x}\tilde{x}}(\bm{x}_{N+1}) =0_{6\times8\times8}
    \end{split}
\end{align}
The numerical implementations are based on \texttt{MATLAB} programming, where the equations of motion are propagated using \texttt{ode113} with relative and absolute tolerances set to $10^{-12}$ (resulting in a local truncation error of $\sim 10^{-7}$ km per integration step). All computations are performed on a personal laptop featuring an 8-core CPU (4 performance cores up to 3.49 GHz, 4 efficiency cores around 2.4 GHz), and 8GB of unified memory. First- and second-order STTs are computed in parallel across the 8 cores. The spacecraft is characterized by a specific impulse of $I_{\mathrm{sp}} = 1950$ s, and an initial mass of $m_{0} = 2000$ kg, while different maximum thrust magnitudes are used, depending on the scenario. Unless stated otherwise, the default algorithm parameters are:
\textit{Trust-region}: $\epsilon_1 = 0.1$, $\epsilon_2 = 0.5$, $\kappa = 0.25$, $\Delta_{\mathrm{max}} = 10^3$, $\Delta_{\mathrm{min}} = 10^{-7}$, and $D = \mathrm{blkdiag}(I_3, 10^2 I_3, 10^{-2})$, where $\mathrm{blkdiag}(\cdot)$ denotes a block-diagonal matrix with the specified blocks on the diagonal;
\textit{Penalty terms}: $\sigma_0 = 10^2$, $\sigma_{\mathrm{max}} = 10^{10}$, $k_{\sigma} = 1.1$, $\delta_{\sigma} = 10^{-5}$;
\textit{Bang-bang smoothing:} initial sharpness $\rho_{1,1} = 1$, step size $k_{\rho} = 0.1$;
\textit{State-path constraint sharpness:} $\epsilon_{\mathrm{SPC}} = 3 \times 10^{-4}$, balancing smoothness and accuracy in the log-sum approximation (see \cref{eq:SmoothApproximationofMax}); a smaller value would cause nonsmoothness near constraint activation, while a larger value would overestimate the $\max\left(\cdot\right)$ function;
\textit{Convergence tolerances:} $\epsilon_{\mathrm{opt}} = 10^{-7}$, $\epsilon_{\mathrm{feas}} = 10^{-7}$, $\rho_{1,\mathrm{min}} = 10^{-2}$.

In the first two numerical examples, we consider scenarios similar to Ref.~\citenum{AzizJGCD} (which employs HDDP) to provide a reference for comparison in terms of fuel consumption, number of DDP stages, and computation time, while the third example is a novel case proposed in this paper.
\subsection{DRO to DRO Transfer: Trade Study Against Indirect Multiple Shooting}
\label{Sec:NumericalResultsDRO}
The first example to demonstrate our PDDP method is a transfer between two Distant Retrograde Orbits (DROs). The initial and final states, time of flight, and maximum thrust are the same as in Sec.VI of Ref.~\citenum{AzizJGCD}. We solve for single- and multi-revolution transfers, increasing the time of flight (TOF) while concurrently decreasing the maximum thrust, without state-path constraints (which is equivalent to set $\beta=0$).

We analyze four distinct scenarios involving trajectories with one, two, five, and ten revolutions. In the first three cases, each trajectory is evenly divided in time into 10 stages, requiring the optimization of $10 \times 7 = 70$ costate variables per scenario. In contrast, the HDDP method used in Ref~.\citenum{AzizJGCD} employs $80$, $160$, and $400$ stages for the same examples, resulting in the optimization of $80 \times 3 = 240$, $160 \times 3 = 480$, and $400 \times 3 = 1200$ variables, respectively. For the ten-revolution trajectory we use $15$ stages, leading to the optimization of $105$ variables, while the approach in Ref.~\citenum{AzizJGCD} utilizes $800$ stages, resulting in $2400$ control variables to optimize. Note that, due to the high sensitivity of the 10-revolution trajectory, a larger smoothing step size of $k_{\rho}=0.5$ was required to ensure convergence to the bang-bang control solution. Using this setup, our PDDP algorithm successfully achieved convergence for each scenario using random initial guess for the costates, demonstrating one of our claims: the analytical law provided by PMP facilitates the design of low-thrust multi-spiral trajectories with a relatively small number of stages. In \cref{tab:DRO2DROMinFuel} we summarize useful results for these single- and multi-revolution transfers, in terms of final mass, number of bang-bang switches, number of PDDP iterations, and computation time. We observe that the cost per iteration for these scenarios is between 0.48 seconds to 0.77 seconds. Despite the need for computing the augmented STTs due to the incorporation of costates in the state vector compared to other DDP methods (see \cref{Sec:StageQuadraticExpansion} for more details), the computation times for these scenarios are similar to those in Ref.~\citenum{AzizJGCD}.

\begin{table}[htbp]
\centering
\begin{tabular}{lccccccc}
\hline Revs &$T_{\mathrm{max}}$, N &$t_{\mathrm{f}}$, days&$m_{\mathrm{f}}$, kg&$\#$ of switch&$\#$ of iter& Runtime, s\\
\hline
  Single&0.25&17.5&1990.86&6&248&120\\
   Two&0.15&35&1993.17&8&618&334\\
   Five&0.05&87.5&1992.54&16&965&740\\
   Ten&0.02&175& 1991.87&20&976&752\\
\hline   
\end{tabular}
\caption{\label{tab:DRO2DROMinFuel} Minimum-fuel DRO to DRO transfer results for different time of flight}
\end{table}

\begin{figure}[hbt!]
\centering \subfigure[\label{fig:DRO2DRO5RevTraj} Optimal trajectory: 5-revolutions transfer]
{\includegraphics[width=0.35\linewidth]{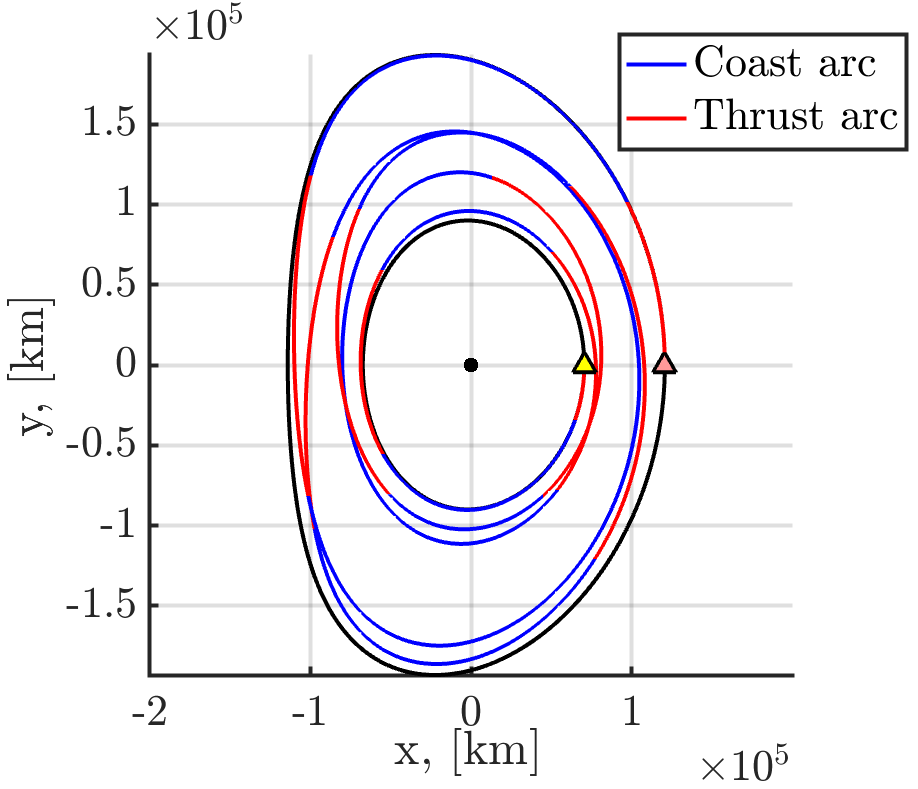}}
\centering \subfigure[\label{fig:DRO2DRO5RevOptimControl} Optimal control: : 5-revolutions transfer]
{\includegraphics[width=0.49\linewidth]{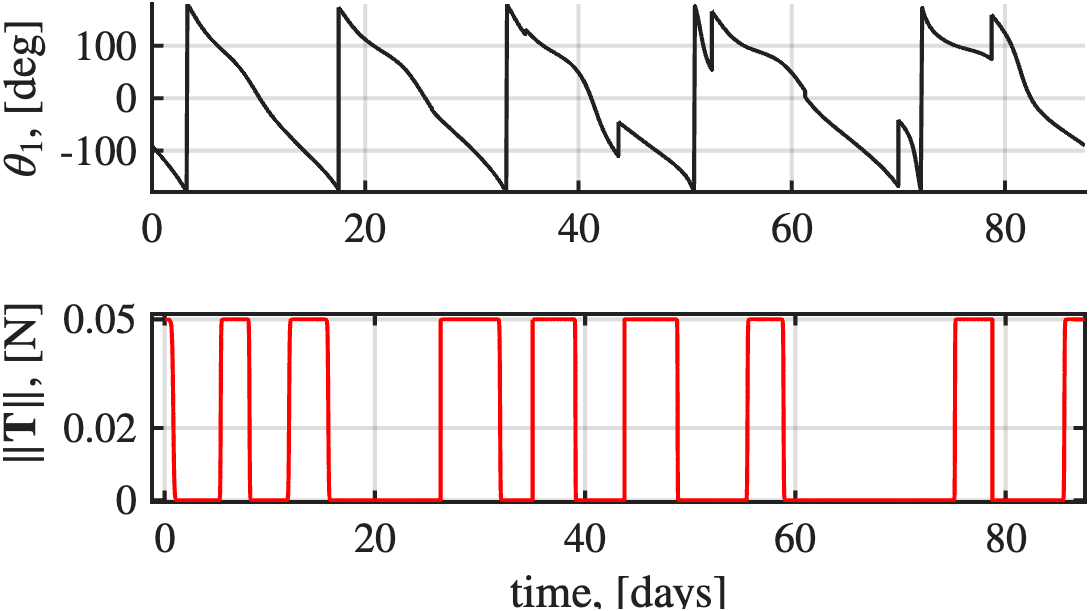}}
\centering \subfigure[\label{fig:DRO2DRO12RevTraj} Optimal trajectory: 10-revolutions transfer]
{\includegraphics[width=0.35\linewidth]{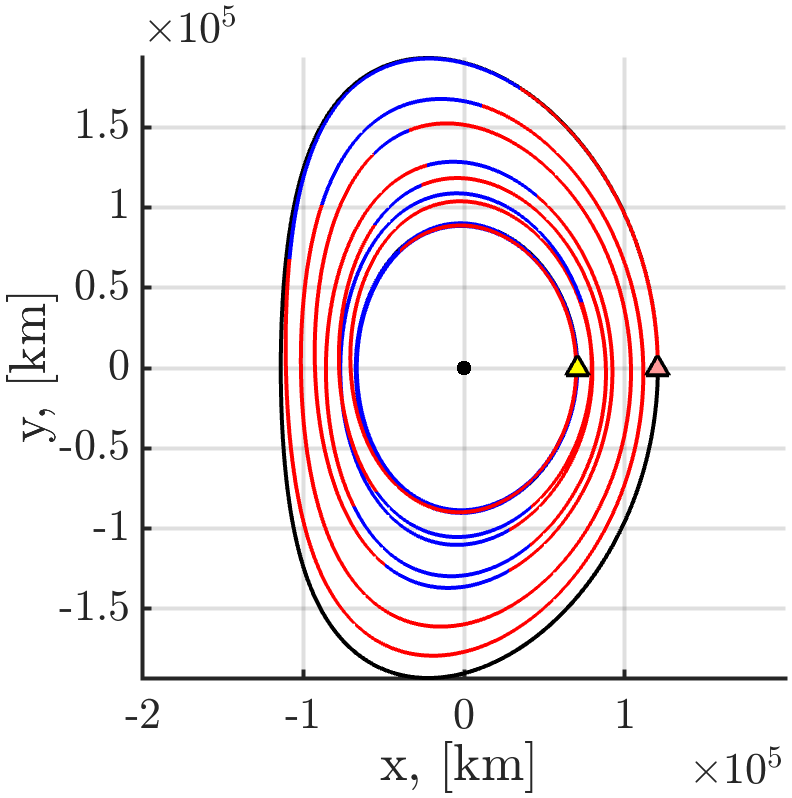}}
\centering \subfigure[\label{fig:DRO2DRO12RevOptimControl} Optimal control: : 10-revolutions transfer]
{\includegraphics[width=0.49\linewidth]{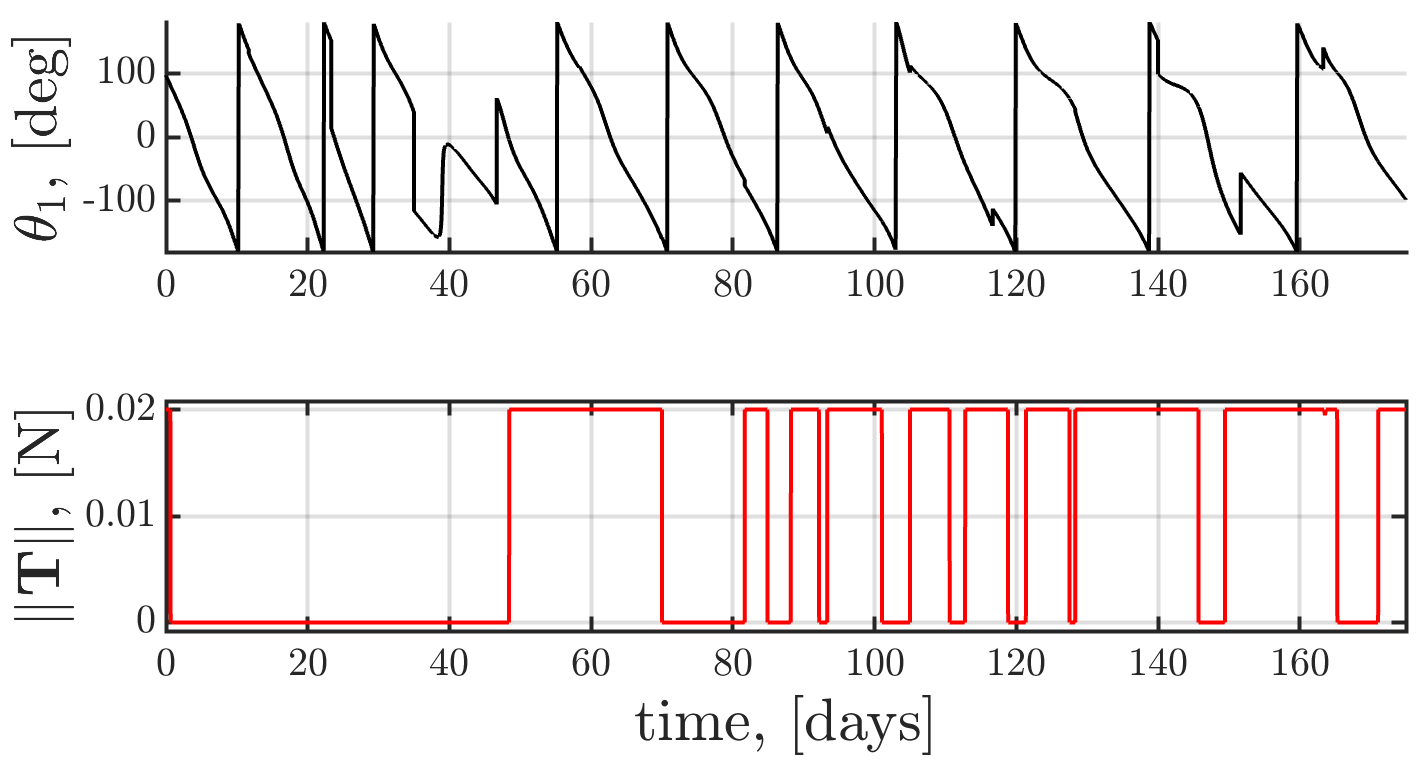}}
\caption{\label{fig:DRO2DROMinFuel} DRO to DRO transfers}
\end{figure}

In \cref{fig:DRO2DRO5RevTraj} we show the five-revolutions minimum-fuel trajectory in the rotating frame where the initial and final states are highlighted in yellow and red, respectively. \cref{fig:DRO2DRO5RevOptimControl} depicts the corresponding optimal control, where $\theta_{1}$, is the longitude of the thrust direction in the rotating frame (the latitude is not plotted, since this is a planar transfer), calculated as: $\theta_{1}=\mathrm{arctan2}(\alpha_{2}/\alpha_{1})$; $\mathrm{arctan2}(.)$ is the four-quadrant inverse tangent function. This solution foresees the same fuel consumption as the one in Ref.~\citenum{AzizJGCD}, i.e. $0.34\%$ of the total mass($\approx 7$ kg). The ten-revolution trajectory in the rotating frame is depicted in \cref{fig:DRO2DRO12RevTraj}, while the corresponding optimal control policy is shown in \cref{fig:DRO2DRO12RevOptimControl}. The solution exhibits 20 switching points for a total fuel consumption of $8.32$ kg (i.e., $0.42\%$ of the total mass). 

Next, to compare the sensitivity of our PDDP algorithm to the initial guess against ``pure'' indirect multiple shooting methods, we analyze the convergence behavior of each approach over 500 samples drawn from a uniform distribution within the interval $[-1000, 1000]$. While previous work proposed costate normalization techniques that confine the initial costate to the interval $[-1, 1]$ \cite{JIANG2012245}, such a technique is not implemented in the present study. Instead, we intentionally adopt a large sampling interval in order to explore the convergence landscape of each method and assess their respective radius of convergence with respect to the initial costate guess. Normalizing the costates to a hyper-sphere would restrict the search space and therefore limit the ability to assess differences in convergence robustness between the approaches. This analysis is performed for the first three scenarios: single-, two-, and five-revolutions, with the smoothing parameter for the bang-bang control set to $\rho_{1}=1$. This value is sufficient for comparison purposes, as solving the smooth-control problem is standard practice and allows continuation to the bang-bang solution \cite{JIANG2012245,SIDHOUMJGCD2024,TAHERI20182470}. The analysis is repeated for different stage counts—using the same number of shooting nodes as PDDP stages—to assess the impact of discretization on sensitivities to poor initial guesses. The indirect multiple shooting method is implemented using \texttt{MATLAB}'s \texttt{fsolve} with analytic gradient to provide more robustness, with a default tolerance of $10^{-6}$. For both approaches, we set the maximum number of iterations to $800$. The results in \cref{tab:DRO2DRORandomGuesses} show that PDDP generally outperforms indirect multiple shooting, with the only exception being the single-revolution transfer with 5 nodes, where the indirect method achieves a slightly higher success rate (75\%) compared to PDDP (72\%). In contrast, while the convergence rate drops below 50\% for the 2- and 5-revolution cases due to the poor quality of random initial guesses, PDDP still outperforms pure indirect shooting. For the two-revolution case with 5 nodes and the five-revolution case with 15 nodes, it achieves success ratios $\sim1.3$ times higher than the indirect method. The most significant improvements are observed in the five-revolution cases with 5 and 10 stages, where PDDP outperforms the indirect method by factors of $\sim1.7$ and $\sim1.8$, respectively. These five-revolution scenarios are particularly noteworthy, as they consistently demonstrate the enhanced robustness of PDDP. Lastly, the trend across the five-revolution tests indicates that PDDP’s robustness improves with an increasing number of stages.

\begin{table}[htbp]
\centering
\begin{tabular}{lccccccc}
\hline & & \multicolumn{2}{c}{\hspace{2cm}Solver (\%) }\\
\cline{3-4} Transfer type &$\#$ of stages (or nodes) &Pure indirect (\textit{fsolve})&PDDP \\
\hline
 \multirow{2}{*}{single-revolution} & 5&75&72\\
   & 10&70.4&75.8\\
   \hline
 \multirow{3}{*}{2-revolution} & 5&36.4&47.2\\
 &10&40.4&40.6\\
   \hline
 \multirow{4}{*}{5-revolution} & 5&20.6&35.2\\
   &10&22.8&41.6\\
   &15&39.8&54.6\\
\hline   
\end{tabular}
\caption{\label{tab:DRO2DRORandomGuesses} Percentage of convergence for 500 random samples}
\end{table}

\subsection{L2 to L1 Halo Orbit Transfer: Minimum-radius Constraint}
\label{Sec:Halo2HaloTransfer}
The next transfer example is an L2 Halo orbit to an L1 Halo orbit with different Jacobi constant, subject to a minimum-radius constraint. Note that, as an alternative to enforcing a path constraint with PDDP, a near-optimal solution could also be obtained by imposing interior boundary conditions on the unconstrained trajectory (e.g., fixing the radius to the minimum value and the radial velocity to zero). The initial and final states, time of flight, and maximum thrust are the same as in Sec.VI of Ref.~\citenum{AzizJGCD}. First, we solve the minimum-fuel transfer \textit{without} imposing any path constraints. The entire trajectory is evenly distributed in time among 3 stages, which involves solving for $3\times7=21$ variables. As a comparison in  Ref.~\citenum{AzizJGCD} the same transfer is discretized into 120 stages, resulting in $120\times3=360$ variables to optimize. Then a constraint on the minimum distance to the Moon is enforced (i.e., $n_{c}=1$), as:
\begin{align}
\label{eq:Example1PathConstraint}
    g(\bm{x}(t),\bm{\lambda}(t)) = r_{2,\mathrm{min}}-r_{2}(t)
\end{align}
where $r_{2}(t) = \sqrt{(x-1+\mu)^{2}+y^{2}+z^{2}}$ is the distance to the Moon, while $r_{2,\mathrm{min}}$ is the user-defined minimum radius. To solve the minimum-fuel trajectory without path constraints, we begin by generating a random initial guess for the costates. The solution of the unconstrained problem serves as the initial guess for the radius-constrained problem with $r_{2,\mathrm{min}}=20,000$ km. Subsequently, we progressively increase the minimum radius constraint to $r_{2,\mathrm{min}}=25,000$ km, and $r_{2,\mathrm{min}}=30,000$ km  using the previous solution as an initial guess for the next step.

\begin{figure}[hbt!]
\centering \subfigure[\label{fig:Halo2HaloTrajNoSPC} Optimal trajectory]
{\includegraphics[width=0.39\linewidth]{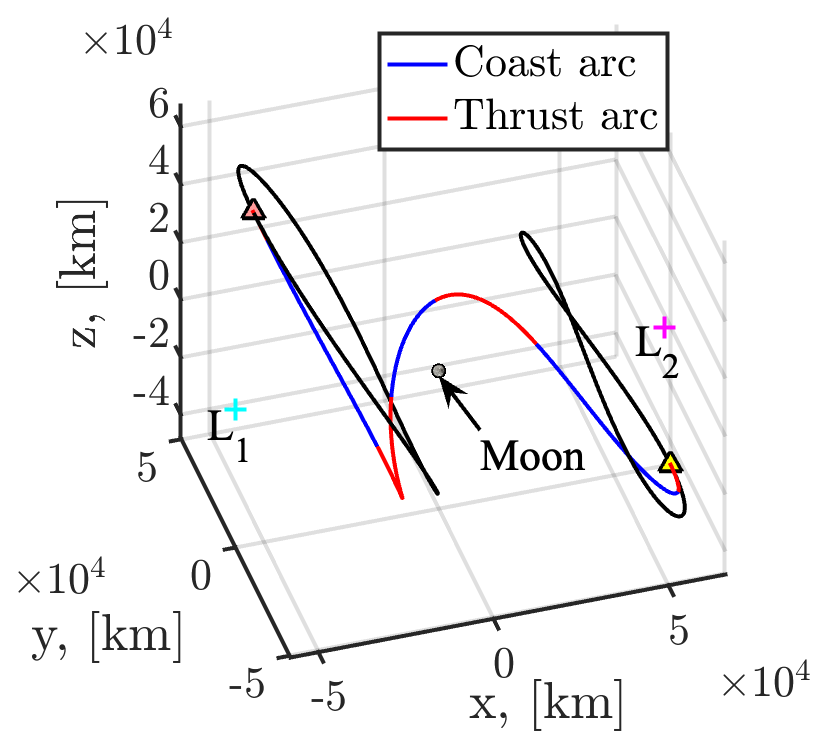}}
\centering \subfigure[Optimal control]
{\includegraphics[width=0.6\linewidth]{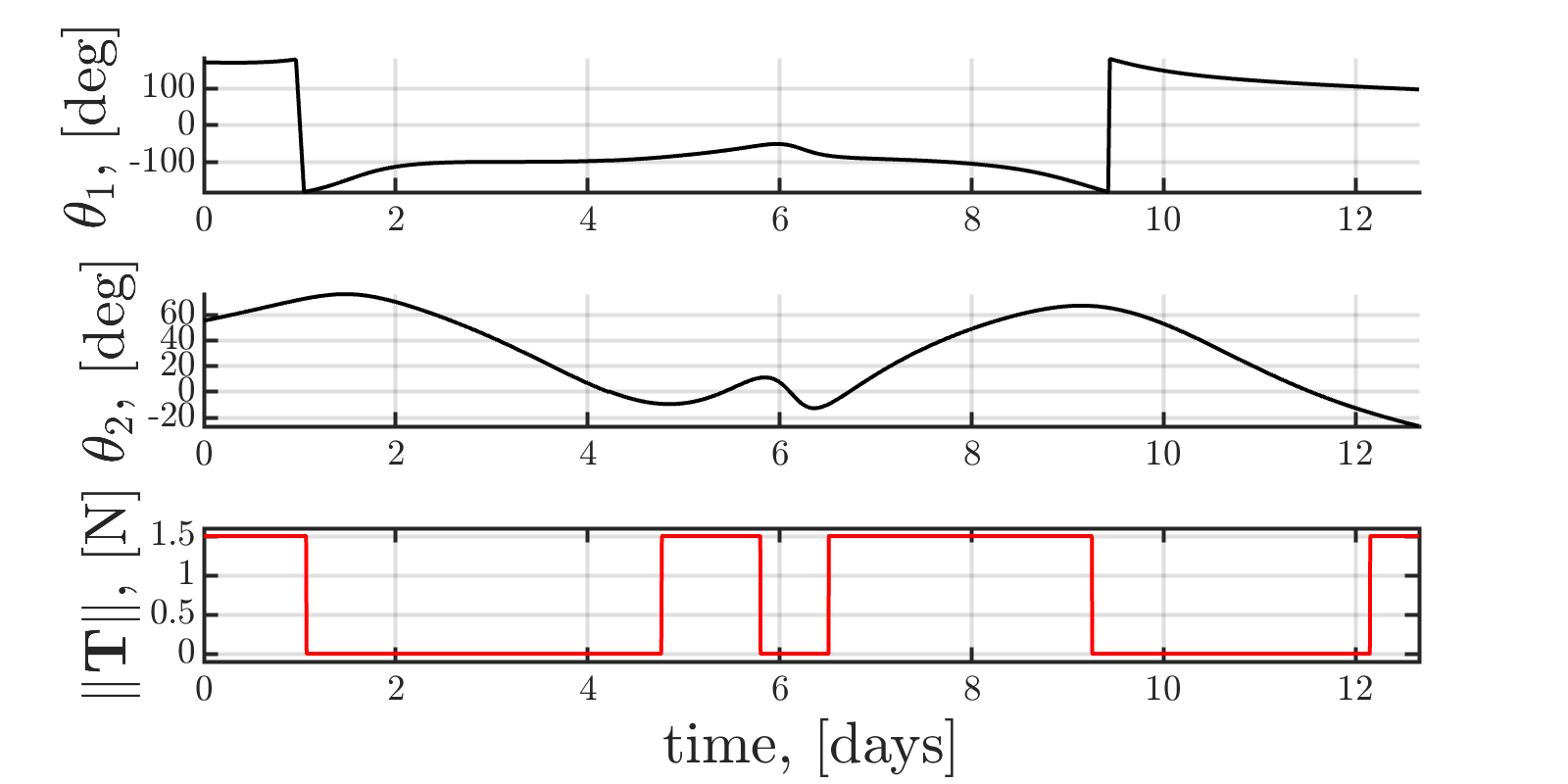}}
\caption{\label{fig:Halo2HaloOptimControl} L2 to L1 Halo orbit minimum-fuel transfer \textit{without} path constraints}
\end{figure}

\begin{figure}[hbt!]
\centering \subfigure[\label{fig:Halo2HaloAllTrajs} Minimum-fuel Halo orbit transfers with increasing radius-constraints]
{\includegraphics[width=0.54\linewidth]{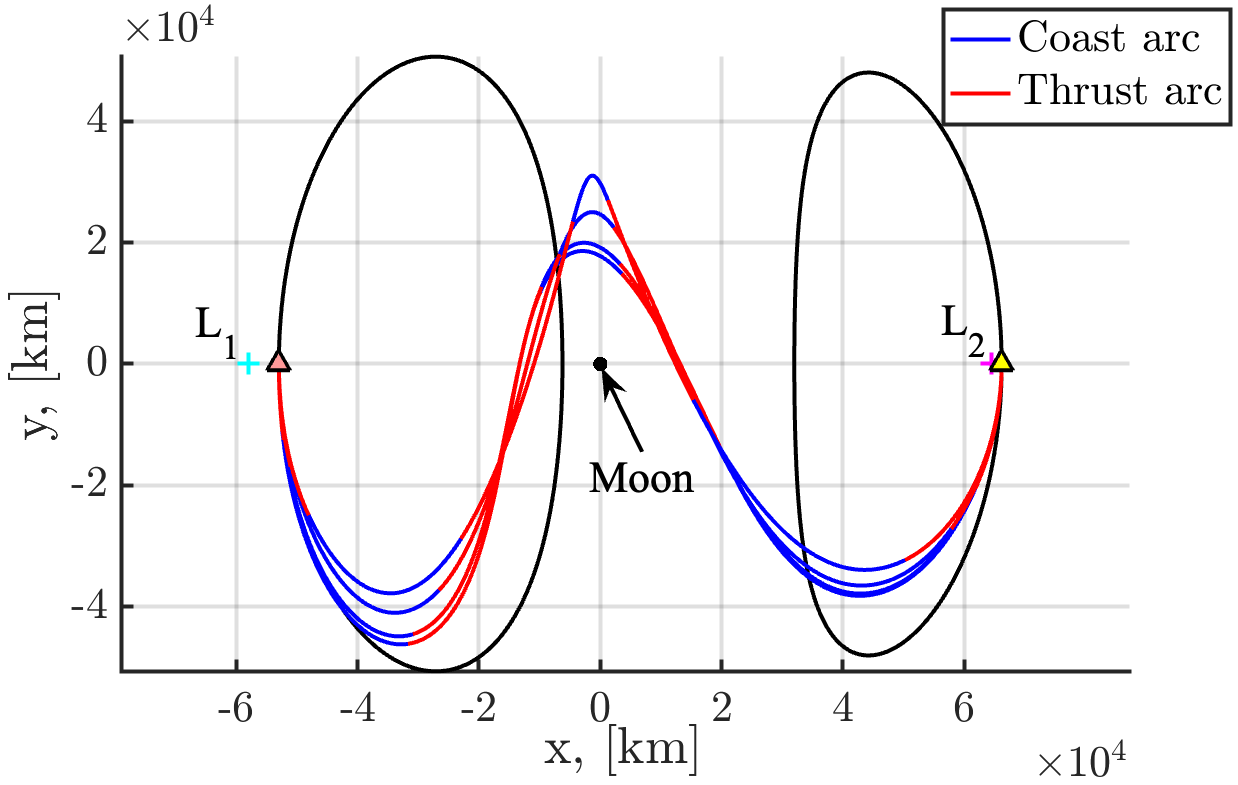}}
\centering \subfigure[\label{fig:Distance2MoonTimeHistory} Distance from the Moon for each radius-constrained transfer]
{\includegraphics[width=0.45\linewidth]{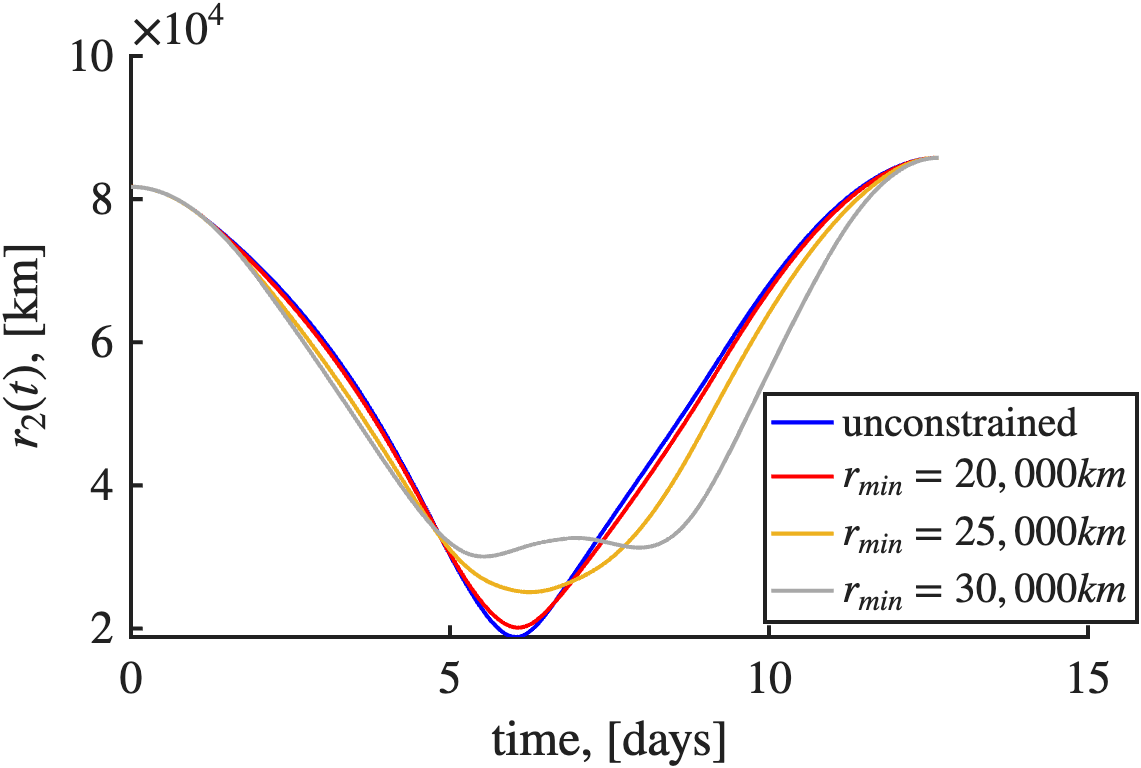}}
\caption{\label{fig:AllTrajHaloHaloOverlaid} L2 to L1 Halo orbit minimum-fuel transfer \textit{with} path constraints}
\end{figure}

The 3D minimum-fuel unconstrained trajectory in the rotating frame along with the corresponding optimal control are depicted in \cref{fig:Halo2HaloOptimControl}, where $\theta_{1}$, $\theta_{2}$ are the longitude and latitude of the thrust direction in the rotating frame, where: $\theta_{2} = \mathrm{arcsin}(\alpha_{3})$. In \cref{fig:Halo2HaloAllTrajs}, we overlay both unconstrained and radius-constrained X-Y minimum-fuel trajectories in the rotating frame. \cref{fig:Distance2MoonTimeHistory} shows the time-history of the distance from the Moon, $r_{2}(t)$, for all four cases. The figure clearly illustrates the effectiveness of our path constraint approach within PDDP, as the minimum distance to the Moon progressively increases with the increasing stringency of the constraint. In \cref{tab:Halo2HaloMinFuel}, we summarize the minimum distance to the Moon, the final mass, the sharpness parameter for the approximation in \cref{eq:SmoothApproximationofMax}, the state-path constraints penalty parameter, the number of iterations, and total computation time. As a reminder, $\beta$ and $\rho_{2}$ are fixed throughout the PDDP iterations. We note that the penalty parameter for the radius-constrained transfer with $r_{2,\mathrm{min}}=30,000$ km is slightly higher than in the other two cases, as this constraint is more stringent and, therefore, more difficult to enforce.
 
 Our solutions correspond to the same local minima as those reported in Ref.~\citenum{AzizJGCD}, albeit with a slightly lower final mass. However, unlike the results in Ref~.\citenum{AzizJGCD}, our approach does not exhibit any path constraint violations, thereby validating the effectiveness of our PDDP method. The solutions are found in a tractable time with an average cost per iteration of $0.6$ seconds, except for the second case, most likely because the newly introduced constraint is initially highly violated, resulting in large penalty gradient near violation and stiffer updates in the quadratic system. Finally, we observe that fuel consumption increases as the minimum-radius constraint becomes more stringent. This outcome aligns with fundamental orbital mechanics principles: without path constraints, the spacecraft takes advantage of close lunar flybys to minimize fuel usage. However, as the minimum radius is increased, additional thrusting arcs on either side of the flyby become necessary.

\begin{table}[htbp]
\centering
\begin{tabular}{lccccccc}
\hline  $r_{2,\mathrm{min}}$, km & $\displaystyle \min\limits_{t} r_2$, km&$m_{\mathrm{f}}$, kg&$\rho_{2}$&$\beta$&$\#$ of iter& Runtime, s\\
\hline  None&18,806&1963.73&$10^{-3}$&$0$&361&165\\   20,000&20,189&1963.62&$10^{-3}$&$10^{3}$&71&82\\ 25,000&25,154&1961.39&$10^{-3}$&$3\times10^{3 }$&
   353&218\\
30,000&30,050&1952.45&$10^{-3}$&$10^{4}$&411&259\\
\hline   
\end{tabular}
\caption{\label{tab:Halo2HaloMinFuel} Minimum-fuel Halo to Halo transfer results for different radius-constraints}
\end{table}

\subsection{9:2 NRHO to 2:1 Retrograde Resonant Orbit Transfer: Eclipse Avoidance Constraint}
\label{Sec:GatewayTransfer}
The final example used to demonstrate our PDDP approach involves a trajectory from the 9:2 Near-Rectilinear Halo Orbit (NRHO) for the Lunar Gateway space station to a 2:1 retrograde resonant orbit, subject to an eclipse-avoidance constraint. Specifically, as discussed in Ref.~\citenum{Gupta2025}, a 2:1 retrograde resonant orbit with a perigee radius equal to that of the geosynchronous orbit (GEO) is identified as a suitable baseline for cislunar observation. A spacecraft in this orbit can repeatedly access the GEO belts while covering the entire cislunar region. This scenario is particularly relevant to NASA’s Artemis program, where a spacecraft could be deployed from the Gateway and travel to this resonant orbit, providing access to the full cislunar domain. The 9:2 NRHO is characterized by a Jacobi constant of $C = 3.0497942$ and a period of $T=6.39$ days, while the 2:1 resonant orbit has a Jacobi constant $C =0.89703010 $ with a period $T=25.96$ days. The significant difference between the Jacobi constants of these two orbits highlights the necessity of low-thrust propulsion to achieve this transfer. The departure is from the apolune of the NRHO, where $\bm{x}_{\mathcal{O},0} = [1.0188262, 0, -0.17979826, 0, -0.09618936, 0]^\top$, and ends on the far-side of the Earth at $y=0$, where $\bm{x}_{\mathcal{O},\mathrm{f}}
= [-1.0280902,0, 0,0,1.4548866,0]$. Note that these initial and final conditions must be adjusted via a differential correction process to obtain truly periodic orbits. The transfer time is fixed and set to $50.85$ days, while the spacecraft maximum thrust is assumed to be $T_{\mathrm{max}} = 3$ N. The time of flight was found by leveraging insights from  multi-body dynamics. We propagated the stable manifolds of the 2:1 resonance and the unstable manifold of the NRHO until they intersected the Poincar\'e section at $x=0.5$.  The closest crossing points were selected, and the corresponding time of flight was used for the transfer. For this transfer, we consider an eclipse-avoidance constraint. We acknowledge that eclipse-aware trajectories can also be optimized by enforcing coasting arcs during eclipse periods (e.g., by multiplying the throttle by an activation function, such as a hyperbolic tangent, where the switching function corresponds to the eclipse constraint). However, strict eclipse-avoidance constraints, as enforced here, remain worthwhile to ensure continuous power availability, which is critical for missions with limited battery capacity or sensitive payloads.

We use the conical shadow model \cite{AzizEclipse}, which assumes spherical shapes of the occulting body (the Earth) and the Sun. For simplicity, the Moon is not modeled as an occulting body—an acceptable assumption, as lunar eclipses predominantly affect transfers within the lunar domain. The conical shadow model defines three parameters, namely the apparent solar radius $\mathrm{ASr}\overset{\underset{\triangle}{}}{=}\mathrm{asin}\left(\frac{R_{\odot}}{\|\bm{r}_{\odot/SC}\|_2}\right)$, the apparent occulting body radius $\mathrm{ABr}\overset{\underset{\triangle}{}}{=}\mathrm{asin}\left(\frac{R_{E}}{\|\bm{r}_{E/SC}\|_2}\right)$, and the apparent distance $\mathrm{AD}\overset{\underset{\triangle}{}}{=}\mathrm{acos}\left(\frac{\bm{r}_{E/SC} \cdot \bm{r}_{\odot/SC}}{\|\bm{r}_{E/SC}\|_2\|\bm{r}_{\odot/SC}\|_2}\right)$, where $R_{\odot}$ and $R_{E}$ are the Sun and Earth physical radii, respectively, while $\bm{r}_{\odot/SC}$ and $\bm{r}_{E/SC}$ are the Sun and Earth relative positions to the spacecraft, respectively. For geometric insights see Fig. 1 in Ref.~\citenum{AzizEclipse}. The spacecraft experiences an eclipse when the following eclipse constraint function is positive:
\begin{align}
\label{Eq:Eclipseconstraint}
    g_{\mathrm{ec}} = \mathrm{ASr}+\mathrm{ABr}-\mathrm{AD}\
\end{align}

\begin{figure}[hbt!]
\centering \subfigure[\label{fig:GatewayTrajNoSPC} X-Y projection]
{\includegraphics[width=0.47\linewidth]{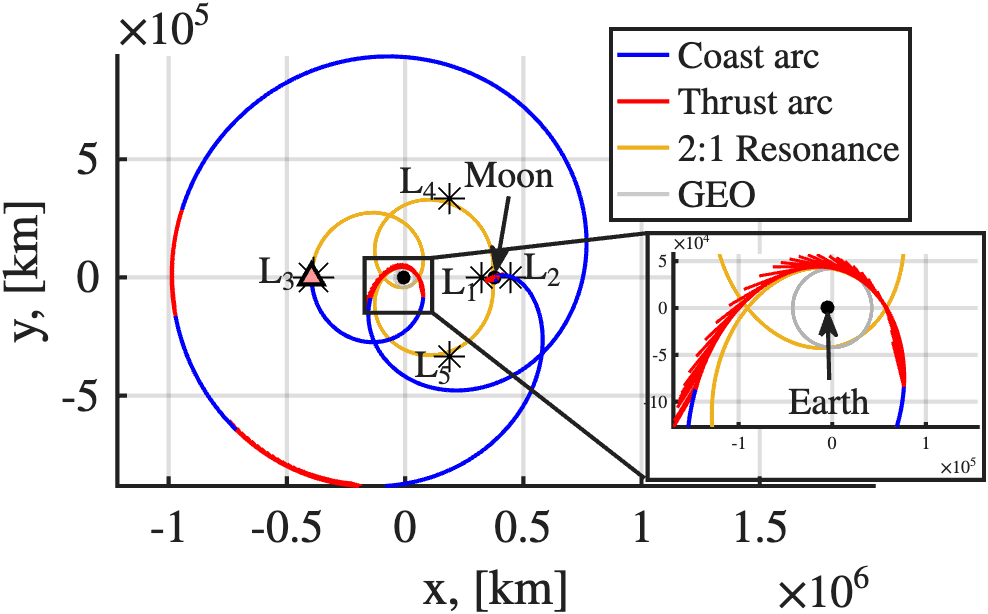}}
\centering \subfigure[\label{fig:GatewayXZTrajNoSPC}X-Z projection (axes not scaled equally)]
{\includegraphics[width=0.52\linewidth]{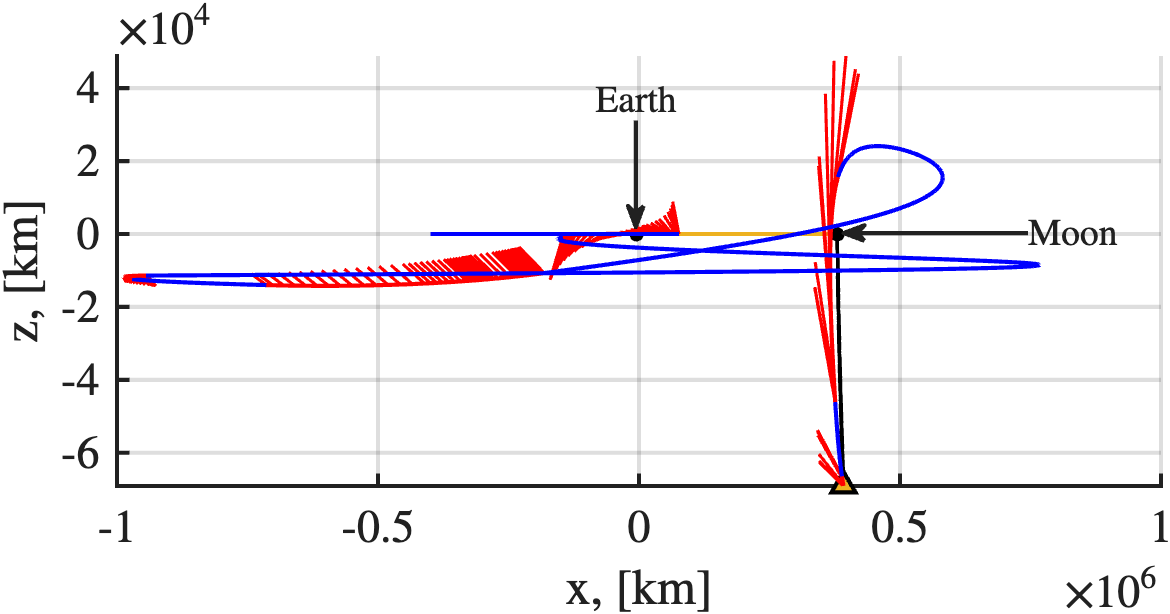}}
\caption{ 9:2 NRHO to 2:1 retrograde resonant orbit: \textit{unconstrained} minimum-fuel transfer}
\end{figure}

\begin{figure}[hbt!]
\centering 
\subfigure[\label{fig:CostateHistory} Costate dynamics: Initial guess vs. Optimal solution]
{\includegraphics[width=0.49\linewidth]{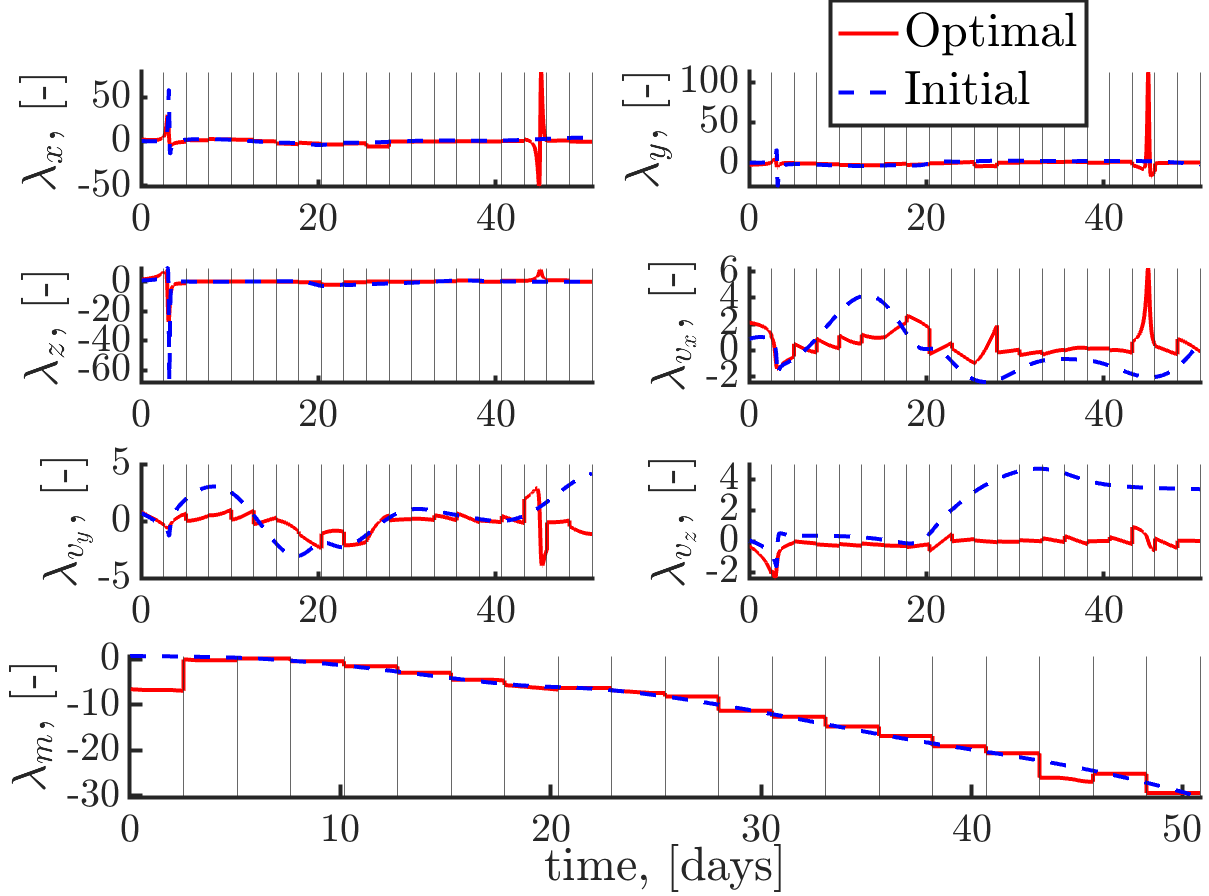}}
\centering 
\subfigure[\label{fig:GatewayOptimControlNoSPC} Optimal control]
{\includegraphics[width=0.49\linewidth]{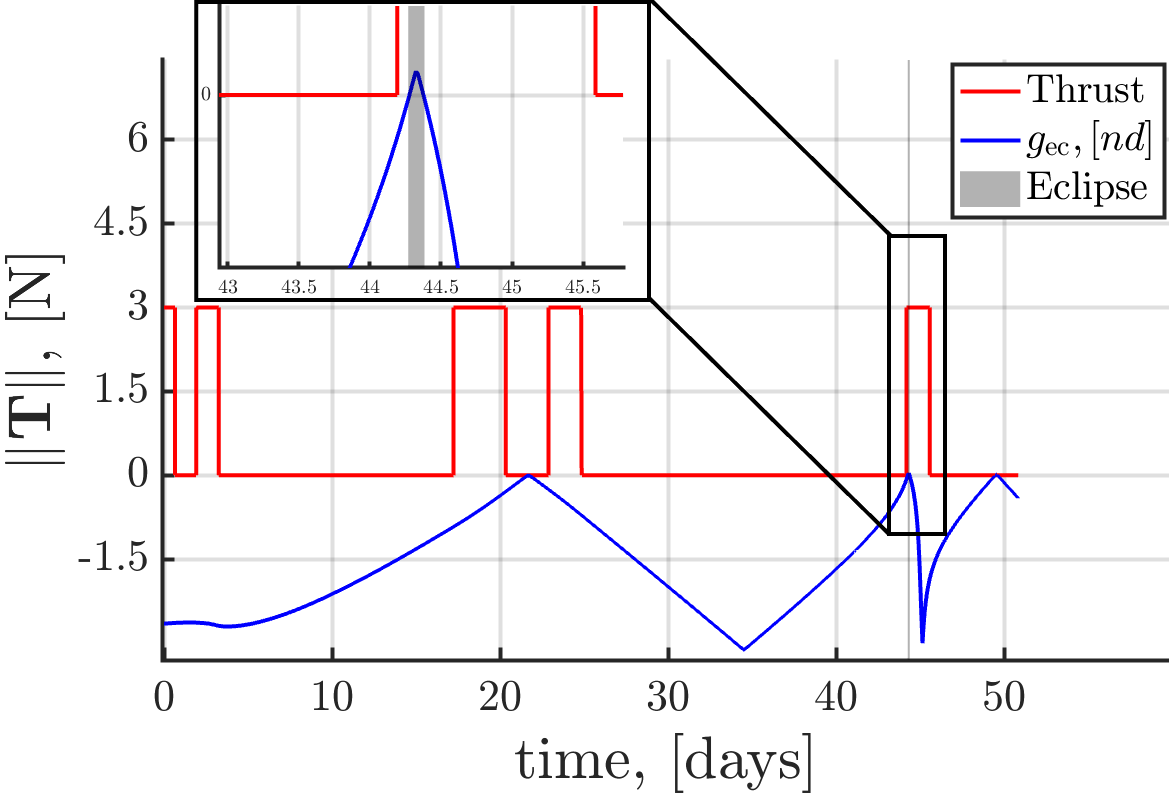}}
\caption{\label{fig:GatewayOptimContAndCostate} 9:2 NRHO to 2:1 retrograde resonant orbit optimal control with shadowing periods}
\end{figure}

Once again, we first solve the unconstrained problem. Due to the complexity of this long-duration scenario we divide the trajectory into $20$ nodes, resulting in $20\times7 = 140$ costate variables to optimize. Using random initial guess for the costates our PDDP algorithm is able to converge to a solution. In \cref{fig:GatewayTrajNoSPC} we show the X-Y projection of the resulting minimum-fuel trajectory, while \cref{fig:GatewayXZTrajNoSPC} depicts the X-Z projection, both in the rotating frame, with the thrust direction highlighted. Note that for the two arcs on the far side of the Earth, the direction of the thrust is difficult to distinguish, as it is nearly tangential to the X-Y trajectory. The spacecraft accomplishes 1 revolution around the Moon before escaping, then performs two revolutions around the Earth, with an Earth flyby in the GEO belt, before capturing the target resonant orbit. This transfer requires $ 5.73\%$ of the total mass (i.e., $114.66$ kg) over the $50.85$ days of flight. In \cref{fig:CostateHistory}, we compare the costate time-history of the optimal solution and the initial guess, with stage locations highlighted by vertical lines. This figure demonstrates the ability of PDDP to effectively adjust the costates along the trajectory towards both feasibility and optimality. Notably, the velocity-costate vector, or primer vector, undergoes the more significant corrections, leading to optimization of both thrust direction and magnitude. Discontinuities in the costates are observed at stage locations; while this could introduce some suboptimality in the PMP-based control law, it is important to note that the quadratic expansion inherent in DDP-based techniques already introduces a degree of suboptimality (see \cref{Sec:StructureofPDDP} for more discussion). \cref{fig:GatewayOptimControlNoSPC} depicts the optimal thrust profile for this transfer with the eclipse function $g_{\mathrm{ec}}$ overlaid. Our solution involves 9 bang-bang switches and predicts an eclipse duration of $\sim165.3$ minutes during the final thrust arc. Eclipse modeling assumes an initial Sun–Earth inertial phasing of $\sim25$ degrees, with Earth in a circular orbit on the X-Y plane around the Sun (inertial frame). This initial phasing was selected because it produces the longest eclipse duration, thereby creating a more challenging, nontrivial constrained optimal control problem. In operation, thrusting during eclipses is infeasible, underscoring the necessity of incorporating eclipse constraints into the optimization framework.

\begin{figure}[hbt!]
\centering \subfigure[\label{fig:GatewayTrajWithSPC} X-Y projection]
{\includegraphics[width=0.47\linewidth]{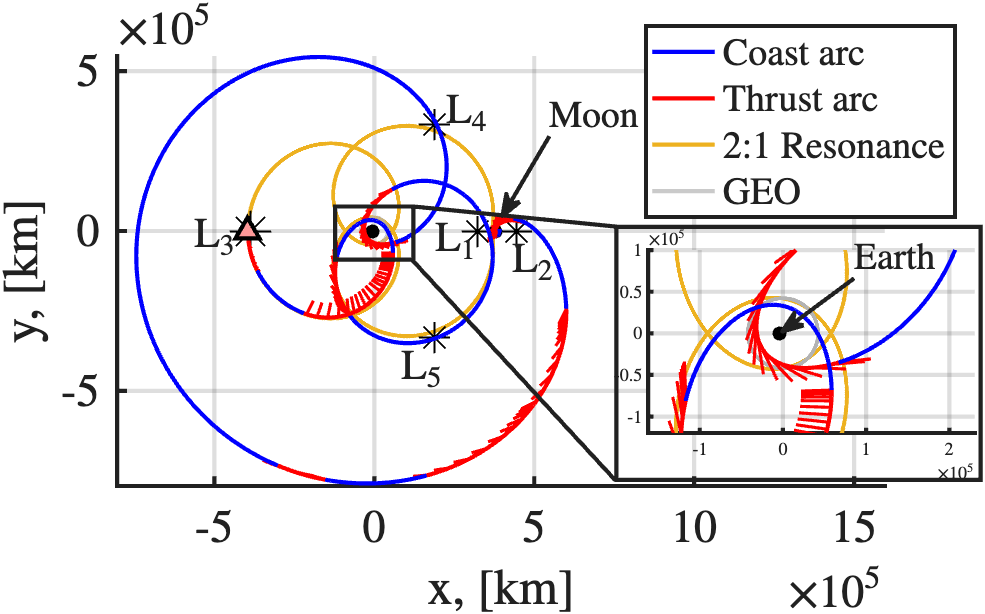}}
\centering \subfigure[\label{fig:GatewayXZTrajWithSPC}X-Z projection (axes not scaled equally)]
{\includegraphics[width=0.52\linewidth]{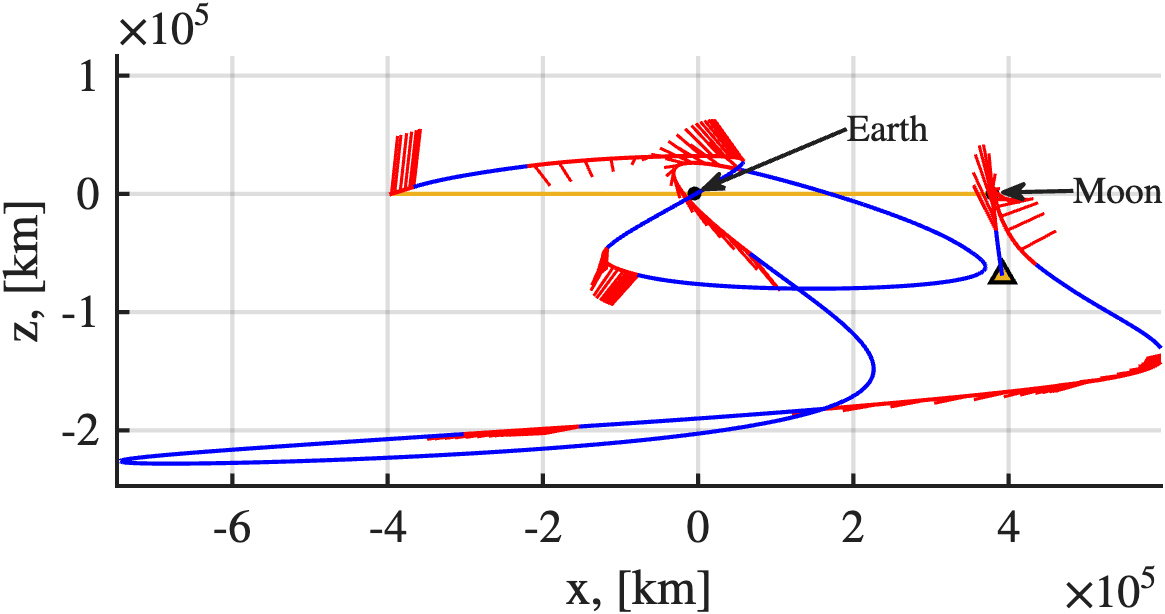}}
\caption{\label{fig:GatewaySolutionWithSPC} 9:2 NRHO to 2:1 retrograde resonant orbit: \textit{eclipse-constrained}  minimum-fuel transfer}
\end{figure}

\begin{figure}[!ht]
\centering
\includegraphics[width=.7\textwidth]{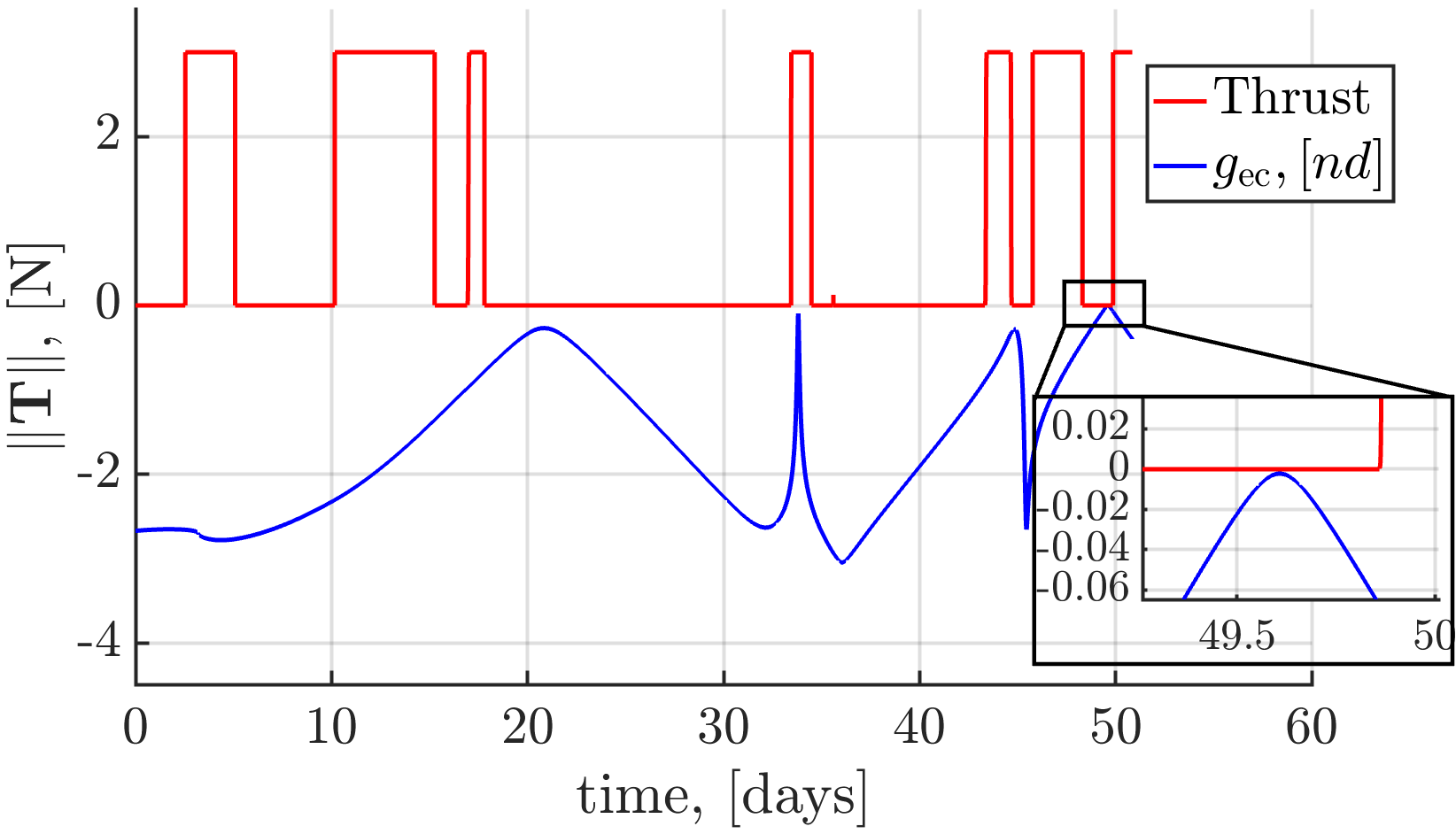}
\caption{9:2 NRHO to 2:1 retrograde resonant orbit optimal control \textit{with} eclipse constraint}
\label{fig:ThrustProfileWithEllipseWithSPC}
\end{figure}

For the eclipse-avoidance constrained optimization, we use a penalty parameter $\beta = 10^3$ and a smoothing parameter for the constraint of $\rho_2 = 2 \times 10^{-4}$. The initial guess is taken from the unconstrained solution with \(\rho_1 = 1\) (smooth control). Because of the high sensitivity of this problem, we use a larger smoothing step size of $k_{\rho}=0.5$ for the bang-bang control. Due to the activation of the eclipse-avoidance constraint, the optimizer converges to a different, eclipse-free, local minimum, as shown in \cref{fig:GatewaySolutionWithSPC}. Notably, the constrained trajectory features two Earth-centered lobes, in contrast to a single lobe in the unconstrained case. The most striking difference appears in the X–Z trajectory shown in \cref{fig:GatewayXZTrajWithSPC}: the constrained trajectory captures the $z = 0$ plane only at the end with a final capture maneuver. Conversely, in the unconstrained case in \cref{fig:GatewayXZTrajNoSPC}, the spacecraft captures the $z = 0$ plane between approximately 44 and 45 days via a thrust arc, followed by an in-plane coast phase. This early capture of the ecliptic plane leads to a solar eclipse, as the spacecraft passes behind the Earth relative to the Sun. In contrast, the constrained solution deliberately avoids capturing the $z = 0$ plane until the end, thereby preventing eclipse by maintaining an out-of-plane trajectory during the critical arc. \cref{fig:ThrustProfileWithEllipseWithSPC} shows the optimal throttle profile with the eclipse constraint overlaid. The solution features 13 bang-bang switches and a consumption of 9.69\% of the total mass ($193.71$ kg), representing an additional $79.7$ kg compared to the unconstrained solution. We note that $g_{\mathrm{ec}}\leq0$ along the entire trajectory, confirming successful eclipse avoidance.

\begin{figure}[hbt!]
\centering \subfigure[\label{fig:OverlayTrajectories} X-Y trajectories]
{\includegraphics[width=0.55\linewidth]{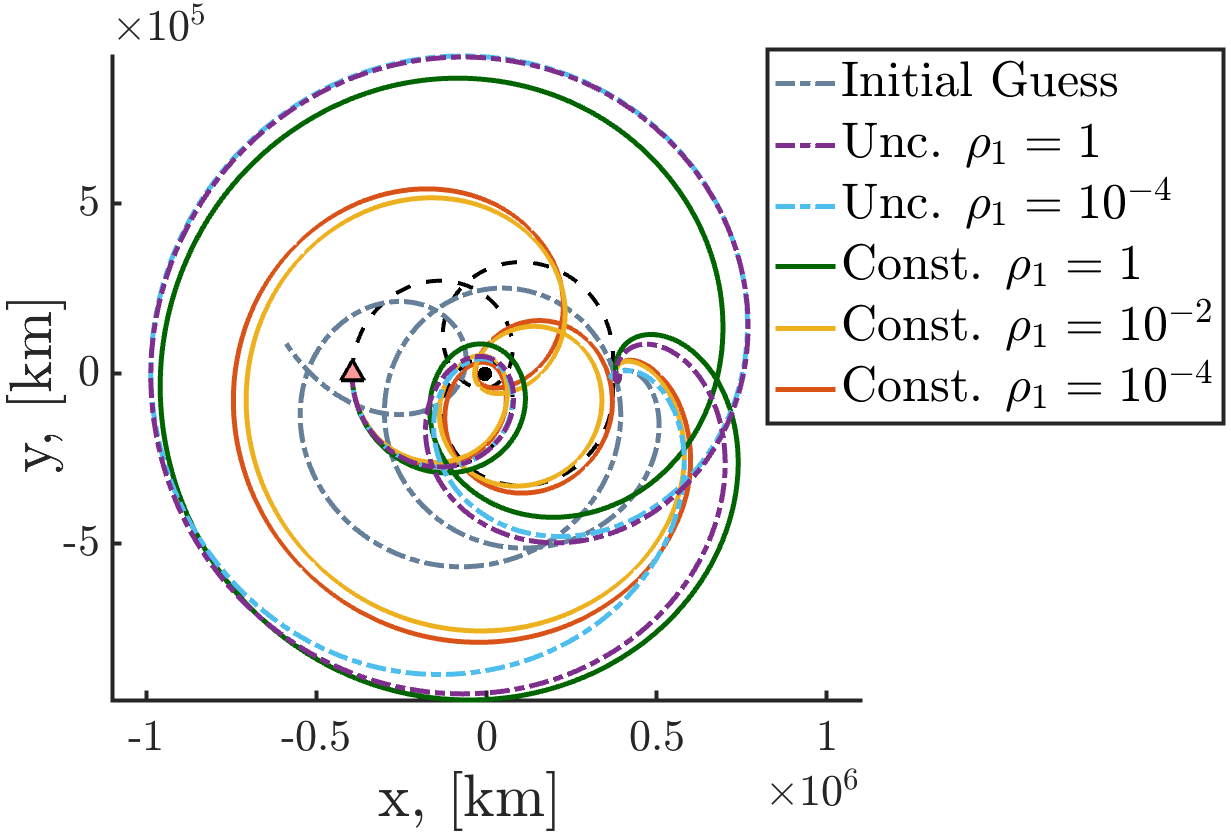}}
\centering \subfigure[\label{fig:OverlayControls}Throttle]
{\includegraphics[width=0.44\linewidth]{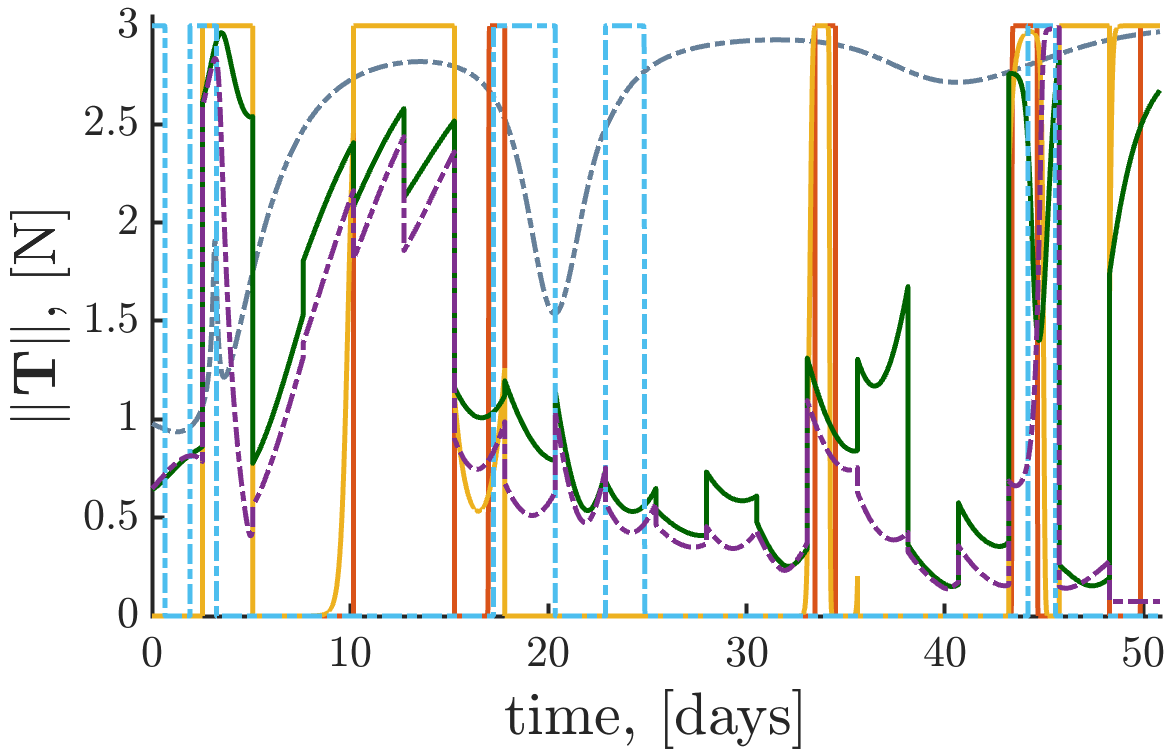}}
\caption{\label{fig:OverlayGateway} 9:2 NRHO to 2:1 Retrograde resonance: Trajectories and controls for varying smoothing parameters}
\end{figure}

Finally, in \cref{fig:OverlayGateway}, we present the trajectories and corresponding control profiles over a decreasing sequence of smoothing parameters, illustrating the solution evolution through successive PDDP iterations. We begin with a random initial guess for the unconstrained problem using a smoothing parameter of $\rho_1 = 1$ (gray curve), which, as expected, does not reach the target. The algorithm then converges to a feasible solution with $\rho_1 = 1$ (purple curve), and this solution is used as the initial guess for the next iteration with a reduced smoothing parameter. This process is repeated, gradually decreasing $\rho_1$ until a bang-bang unconstrained solution is obtained at $\rho_1 = 10^{-4}$ (light blue curve). Next, we consider the eclipse-constrained problem. The smooth unconstrained solution (with $\rho_1 = 1$) is used as the initial guess, and the algorithm converges to an eclipse-free solution at $\rho_1 = 1$ (green curve). As before, the smoothing parameter is progressively decreased, ultimately yielding a bang-bang, eclipse-free trajectory at $\rho_1 = 10^{-4}$ (red-orange curve). In \cref{fig:OverlayTrajectories}, we observe that as the solution transitions toward a bang-bang, eclipse-free control structure, the escape spiral from the Moon becomes progressively tighter. This figure also highlights the gradual transition from a single-lobed Earth-centered trajectory to a two-lobed structure. The transition from smooth control to bang-bang solution, as well as the difference in control profiles between constrained and unconstrained solutions can be appreciated in \cref{fig:OverlayControls}.

\subsection{Discussion}
The three numerical examples highlight the key advantages of incorporating PMP into the DDP formulation. First, the trade study against `pure' indirect multiple shooting indicates that PDDP achieves better convergence robustness, with higher convergence ratios over large initialization intervals in most of the cases. Second, we demonstrated the effectiveness of our approach to handle continuous-time state-path constraints through a minimum-radius constraint in a Halo-to-Halo transfer and an eclipse avoidance constraint in a NRHO to 2:1 resonant orbit transfer. 
Third, the first two examples allowed comparison with state-of-the-art DDP solvers, confirming that PDDP computes optimal solutions in comparable time and efficiently generates long-duration, multi-spiral low-thrust trajectories with fewer optimization variables. Although the present paper focuses on trajectories with several revolutions, the framework is fully applicable to cases with many more revolutions (e.g., hundreds near Earth), though a larger number of PDDP stages would likely be needed to ensure convergence.

As we discuss room for improvements, a promising alternative to pure penalty methods for enforcing continuous-time state-path constraints is the use of augmented Lagrangian techniques—similar to those used for terminal equality constraints—which can alleviate numerical issues such as ill-conditioning, increased nonlinearity, and slow convergence associated with large penalty weights \cite{Fletcher2000}. In addition implementing automated or adaptive tuning strategies would still further enhance the algorithm robustness. Nonetheless, we provide the following recommendations for tuning these parameters: penalty parameters should not become too large to avoid ill-conditioning, so $\beta$ and $\sigma_{\mathrm{max}}$ should be adjusted accordingly. The update factor $k_{\sigma}$ in \cref{eq:PenaltyParameterUpdate} balances conditioning and convergence speed: too large may overly increase fuel consumption to aggressively reduce constraints violation, while too small slows PDDP convergence. The same applies to $\delta_{\sigma}$, which should not increase $\sigma$ unnecessarily once constraints violation are significantly reduced (i.e., avoid $\delta_{\sigma}$ too large). In terms of updating the bang-bang smoothing parameter with $k_{\rho}$ in \cref{eq:SmoothingParameterUpdate}, smaller values reduce the number of steps to capture a bang-bang profile, however can cause large control changes, significantly increasing constraints violation. Larger values, on the other hand, ensure smoother transitions at the cost of convergence speed.Finally, note that the computational cost for computing the STTs could be further reduced by leveraging sparsity. 
\section{Conclusion}
\label{Sec:Conclusion}
This paper presents a Differential Dynamic Programming (DDP) framework that incorporates Pontryagin’s Minimum Principle (PMP) to parameterize the control, for constrained low-thrust trajectory optimization, termed as Pontryagin-Bellman Differential Dynamic Programming (PDDP). 

The presented PDDP algorithm offers three key advantages over existing methods. First, the capability to handle general continuous-time constraints by integrating violations along the trajectory and penalizing them quadratically—a feature that pure indirect methods can hardly accommodate. This approach is general, does not rely on primer vector, and is applicable to any DDP formulation. Second, the improved robustness to poor initial guesses by using DDP to optimize costate variables, compared to traditional indirect shooting approaches. Third, the use of an analytical control law derived from PMP enables efficient generation of multi-spiral low-thrust trajectories with fewer optimization variables than existing DDP-based solvers.

These features are demonstrated via three numerical examples in the Earth-Moon system. The flexibility of PDDP in handling various state-path constraints is demonstrated through minimum-radius and eclipse-avoidance trajectories. A trade study against indirect multiple shooting methods shows that PDDP achieves higher convergence rates overall; notably, for a multi-revolution transfer, it outperforms the indirect multiple shooting method, highlighting improved robustness to poor initial guesses. While the focus is on low-thrust trajectory optimization, the general formulation makes PDDP applicable to any control-affine optimal control problem.
\appendix
\section{Proof of Lemma 1}
\label{Sec:ProofLemma1}

\textit{Proof:} The first property (1) is shown first. Provided that the state-path constraints $g(\cdot)$ is twice continuously differentiable with respect to $\bm{x}(t)$ and $\bm{\lambda}(t)$, it follows from \cref{eq:SmoothApproximationofMax} that $\tilde{l}(\cdot)$ is twice continuously differentiable with respect to $\bm{x}(t)$ and $\bm{\lambda}(t)$, since both $\log(\cdot)$ and $\exp(\cdot)$ are also twice continuously differentiable (these functions are even infinitely continuously differentiable).  The second property (2) is verified by showing that $\tilde{l}(\bm{x}(t),\bm{\lambda}(t))\geq l(\bm{x}(t),\bm{\lambda}(t))$. First since $\log(\cdot)$ is strictly increasing  $\log{\left(1+\exp{\left(\frac{g(\bm{x}(t),\bm{\lambda}(t))}{\rho_{2}}\right)}\right)}\rho_{2} \geq \log{\left(\exp{\left(\frac{g(\bm{x}(t),\bm{\lambda}(t))}{\rho_{2}}\right)}\right)}\rho_{2} = g(\bm{x}(t),\bm{\lambda}(t))$. On the other hand since $\exp(a)>0$ $\forall a$, $\log{\left(1+\exp{\left(\frac{g(\bm{x}(t),\bm{\lambda}(t))}{\rho_{2}}\right)}\right)}\rho_{2} \geq \log{(1)}\rho_{2}=0$. These two inequalities guarantee that $\tilde{l}(\bm{x}(t),\bm{\lambda}(t))\geq \max\left[0,g(\bm{x}(t),\bm{\lambda}(t))\right]=l(\bm{x}(t),\bm{\lambda}(t))$. The third property (3) is verified by showing the following three equations: a) $ \displaystyle \lim_{ \rho_{2}\to 0^{+}}\tilde{l}(\bm{x}(t),\bm{\lambda}(t)) =0$ if $g(\bm{x}(t),\bm{\lambda}(t))<0$; b) $ \displaystyle \lim_{ \rho_{2}\to 0^{+}}\tilde{l}(\bm{x}(t),\bm{\lambda}(t)) =0$ if $g(\bm{x}(t),\bm{\lambda}(t))=0$; c)$ \displaystyle \lim_{ \rho_{2}\to 0^{+}}\tilde{l}(\bm{x}(t),\bm{\lambda}(t)) =g(\bm{x}(t),\bm{\lambda}(t))$ if $g(\bm{x}(t),\bm{\lambda}(t))>0$. Equations (a) is shown by noting that $\displaystyle \lim_{ \rho_{2}\to 0^{+}}\frac{g(\bm{x}(t),\bm{\lambda}(t))}{\rho_{2}} =-\infty$ if $g(\bm{x}(t),\bm{\lambda}(t))<0$, that $\displaystyle \lim_{ a\to -\infty}\exp(a) =0$, and that $\log(\cdot)$ is continuous. Equation (b) is shown by noting that $\displaystyle \lim_{ \rho_{2}\to 0^{+}}\log(2)\rho_{2}=0$. Equation (c) is shown by noting that $\displaystyle \lim_{ \rho_{2}\to 0^{+}}\frac{g(\bm{x}(t),\bm{\lambda}(t))}{\rho_{2}} =+\infty$ if $g(\bm{x}(t),\bm{\lambda}(t))>0$, and that:  
\begin{align}
\begin{split}
        \log{\left(1+\exp{\left(\frac{g(\bm{x}(t),\bm{\lambda}(t))}{\rho_{2}}\right)}\right)}\rho_{2}
    &=g(\bm{x}(t),\bm{\lambda}(t))+\log{\left(1+\exp{\left(-\frac{g(\bm{x}(t),\bm{\lambda}(t))}{\rho_{2}}\right)}\right)}\rho_{2}\\
    &=g(\bm{x}(t),\bm{\lambda}(t))+\mathcal{O}\left(\rho_{2}\right)
    \end{split}
\end{align}
This completes the proof. \hfill $\square$
\section*{Acknowledgments}
This material is based upon work supported in part by the Air Force Office of Scientific Research under award number FA9550-23-1-0512.

%%       \citet{<label>} ==> Lamport [21]
%%       \citep{<label>} ==> [21]
%%




\end{document}